\documentclass[12pt]{article}

\usepackage{amsmath}        
\usepackage{latexsym}

\newtheorem{teor}{Theorem}[section]

\newtheorem{cor}{Corollary}[section]
\newtheorem{obs}{Remark}[section]
\newtheorem{defin}{Definition}[section]

\newtheorem{prov}{Proof of Theorem}[section]

\newtheorem{algor}{Algorithm}[section]

\newfont{\Mb}{msbm10}

\begin{document}
\setcounter{equation}{0}
\setcounter{figure}{0}
\setcounter{table}{0}

\hspace\parindent
\thispagestyle{empty}

\bigskip
\bigskip
\bigskip
\begin{center}
{\LARGE \bf Finding nonlocal Lie symmetries}
\end{center}
\begin{center}
{\LARGE \bf  algorithmically}
\end{center}

\bigskip

\begin{center}
{\large
$^a$L.G.S. Duarte, $^a$L.A.C.P. da Mota and $^a$A.F. Rocha \footnote{E-mails: lgsduarte@gmail.com, lacpdamota@gmail.com and alexfisica1976@yahoo.com.br}
}

\end{center}

\bigskip
\centerline{\it $^a$ Universidade do Estado do Rio de Janeiro,}
\centerline{\it Instituto de F\'{\i}sica, Depto. de F\'{\i}sica Te\'orica,}
\centerline{\it 20559-900 Rio de Janeiro -- RJ, Brazil}

\bigskip\bigskip
\bigskip
\bigskip

\abstract{Here we present a new approach to compute symmetries of rational second order ordinary differential equations (rational 2ODEs). This method can compute Lie symmetries (point symmetries, dynamical symmetries and non-local symmetries) algorithmically. The procedure is based on an idea arising from the formal equivalence between the total derivative operator and the vector field associated with the 2ODE over its solutions (Cartan vector field). Basically, from the formal representation of a Lie symmetry it is possible to extract information that allows to use this symmetry practically (in the 2ODE integration process) even in cases where the formal operation cannot be performed, i.e., in cases where the symmetry is nonlocal. Furthermore, when the 2ODE in question depends on parameters, the procedure allows an analysis that determines the regions of the parameter space in which the integrable cases are located.}

\bigskip
\bigskip
\bigskip
\bigskip
\bigskip
\bigskip


{\bf MSC class: 34-xx}

{\bf ACM class: G.1.7 }

\newpage

\section{Introduction}
\label{intro}

The Lie symmetry method is, probably, the most powerful method to search for first integrals of ordinary differential equations (ODEs). The great majority of integration techniques can be viewed as particular cases of a general method of integration based on the extended group of symmetries admitted by the differential equation (DE). The main obstacle to the application of the Lie method is, ironically, the very computation of symmetries because, up to date, there is no systematic way to find the symmetries of an ODE in the general case. Lie himself has created some procedures for particular cases: for example, in the case of second order ODEs (2ODEs) he considered their invariance through point symmetries, i.e., a symmetry that depends only on the dependent and independent variables. This assumption allows the partial differential equation (PDE) which determines the symmetry condition (Determining equation) to be `separated' in powers of the derivative, so that this PDE would result in an over determined system of PDEs that (in principle) could be solved. But, in the case of an ODE presenting only dynamical or non-local symmetries this strategy cannot be applied. So, many approaches were developed to overcome this difficulty: in \cite{Noscpc1998} Cheb-Terrab {\it et al} developed some heuristics that can find dynamical symmetries of second order ODEs; P. J. Olver introduced the concept of exponential vector field (see \cite{Olv}, p. 185); B. Abraham-Shrauner, A. Guo, K.S. Govinder, P.G.L. Leach, F.M. Mahomed, A.A. Adam, M. L. Gandarias, M. S. Bruzón, M. Senthilvelan and others worked with the concept of hidden and non-local symmetries \cite{AbrGuo,AbrGuo2,AbrGovLea,Abr,GovLea,AdaMah,GanBru,GanBruSen}. C. Muriel and J.L. Romero have developed the concept of $\lambda$-symmetry \cite{MurRom,MurRom2}. See also \cite{MurRom3,MurRom4,CicGaeWal,CicGaeMor,Nuc,Nuc2}. However, despite all these efforts, we still not have a general algorithm to find the symmetries even in more specific cases. 

In this work we show a semi algorithmic procedure to compute the symmetries (whether point, dynamical or non-local) of a rational 2ODE. The main advantage of our method is its algorithmic nature and, on the other side, the principal disadvantage is that it applies only to rational 2ODEs that present a Liouvillian first integral such that their derivatives are of the form $\,{\rm e}^{A/B}\prod p_i^{n_i}$, where $A$, $B$ and $p_i$ are polynomials in $(x,y,y')$ and the $n_i$ are constants. However, we hope to show throughout the paper that this restriction is not as strong as it seems to be at first.

This paper is organised as follows: 

In section \ref{bafac}, we list some basic facts about solving/reducing ODEs using Lie symmetry method.

In section \ref{themet}, we present some results that allow for the construction of a procedure. Then, we introduce the steps of the algorithm and an example in order to clearer the application of the method. 

In section \ref{speci}, we show some cases where the algorithm can be quite improved.

In section \ref{perfor}, we apply the method to the non linear rational 2ODEs of Kamke's book and to some 2ODEs presenting non local symmetries and discuss the algorithm's performance.

In section \ref{applic}, we apply the method to some 2ODEs that describe physical systems.

Finally, we present our conclusions and directions to further our work.

\section{Some basic results about Lie symmetry method for 2ODEs}
\label{bafac}

In this section we present some basic results that will be used in the following discussion.

\subsection{If we have only one point symmetry}
\label{ops}

Consider a rational 2ODE 
\begin{equation}
\label{r2ode}
y'' = {\frac{M(x,y,y')}{N(x,y,y')}} = \phi(x,y,y'),
\end{equation}
where $M$ and $N$ are co prime polynomials in $(x,y,y')$\footnote{$f'$ denotes $\frac{df}{dx}$.}, that presents a global analytical first integral $I(x,y,y')$\footnote{In other words, an integrable 2ODE (in the analytical sense).}. 

The steps of the Lie method are:

\begin{enumerate}
\item Find a continuous group $G$ that is a symmetry group\footnote{$G$ is a group of transformations that does not change the format of the 2ODE.} of the 2ODE (\ref{r2ode}). The conditions for $G$ to be a a symmetry group of the 2ODE (\ref{r2ode}) are usually expressed in one of the following ways:
\begin{enumerate}
\item First way: 
\begin{equation}
\label{condsym1}
X^{(2)}\left[y''-\phi(x,y,y')\right] = 0, \,\,\,\, {\rm mod}\left(y''-\phi(x,y,y')=0\right),
\end{equation}
where $X$ is the symmetry vector field (commonly called group generator) and $X^{(2)}$ is its second extension:
\begin{eqnarray}
X &\equiv& \xi(x,y)\,\frac{\partial}{\partial x}+\eta(x,y)\,\frac{\partial}{\partial y}, \nonumber \\  [3mm]
X^{(2)} &=& \xi(x,y)\,\frac{\partial}{\partial x}+\eta(x,y)\,\frac{\partial}{\partial y}+\eta^{(1)}(x,y,y')\,\frac{\partial}{\partial y'}+\eta^{(2)}(x,y,y',y'')\,\frac{\partial}{\partial y''}, \nonumber
\end{eqnarray}

\noindent
and $\eta^{(i)}=\left[ \displaystyle{{\frac {d }{d x}}
\eta^{(i-1)}} - y^{(i)}\displaystyle{{\frac {d }{d x}}\xi} \right]$.

\begin{obs} \

\noindent
(i) $\eta^{(i)}$ is linear in $y^{(i)}$;

\noindent
(ii) $\eta^{(i)}$ is a polynomial in $(y', y'',\ldots,y^{(i)})$ with coefficients linear and homogeneous in $\xi(x,y), \eta(x,y)$ and their partial derivatives up to the $i^{th}$ order.
\end{obs}

\item Second way:
\begin{equation}
\label{condsym2}
\left[D_x,X^{(1)}\right] = D_x\left[\xi\right] D_x,
\end{equation}
where $X^{(1)}$ is the first extension of $X$
\begin{equation}
X^{(1)} = \xi(x,y)\,\frac{\partial}{\partial x}+\eta(x,y)\,\frac{\partial}{\partial y}+\eta^{(1)}(x,y,y')\,\frac{\partial}{\partial y'},
\end{equation}
and $D_x$ is the total derivative with respect to the independent variable ($x$) over the solutions of the 2ODE (\ref{r2ode}), i.e.,
\begin{equation}
D_x \equiv \frac{d}{dx}{\bigg|}_{y''=\phi(x,y,y')}.
\end{equation}
Since $\eta^{(1)}$ depends on $(x,y,y')$, then $\displaystyle{\frac{d}{dx} = \frac{\partial}{\partial x}+y'\,\frac{\partial}{\partial y}+y''\,\frac{\partial}{\partial y'}}$ and so, $\displaystyle{D_x =\frac{\partial}{\partial x}+y'\,\frac{\partial}{\partial y}+\phi(x,y,y')\,\frac{\partial}{\partial y'}}$.

\end{enumerate}

\vspace{3mm}
\item Reduce the order of the 2ODE using the symmetry:
\begin{enumerate}
\item Using {\em canonical coordinates}:
In the coordinates $(r,s)$ defined by
\begin{eqnarray}
\xi(x,y)\,\frac{\partial r}{\partial x}+\eta(x,y)\,\frac{\partial r}{\partial y}=0, \nonumber \\ [3mm]
\xi(x,y)\,\frac{\partial s}{\partial x}+\eta(x,y)\,\frac{\partial s}{\partial y}=1, \nonumber
\end{eqnarray}
the 2ODE will not depend on $r$ explicitly: $\displaystyle{\frac{d^2 s}{dr^2}={\varphi}\left(s,\frac{ds}{dr}\right)}$.

\item Using {\em differential invariants}:
Every 2ODE that admits the group $G$ can be written as
\begin{equation}
\frac{dI^{(1)}}{dI^{(0)}} = F(I^{(0)},I^{(1)}),
\end{equation}
where ${I^{(0)}}$ and ${I^{(1)}}$ are the zeroth and first order invariants of the group $G$, respectively. So, in the coordinates $(u,v)$ defined by
\begin{eqnarray}
u=I^{(0)}, \nonumber \\ [3mm]
v=I^{(1)}, \nonumber
\end{eqnarray}
the 2ODE will be written as $\displaystyle{\frac{dv}{du}=\psi\left(u,v\right)}$.

\end{enumerate}
\end{enumerate}

\medskip

\subsection{If we have two point symmetries}
\label{tps}

\begin{enumerate}
\item Find two symmetry transformations: In this case we have two (independent) symmetry groups whose generators are given by
\begin{eqnarray}
X_{1} & = & \xi_{1}{\frac{\partial }{\partial x}}+
\eta_{1}{\frac{\partial }{\partial y}}+
\eta_{1}^{(1)}{\frac{\partial }{\partial y'}}\:, \\
X_{2} & = & \xi_{2}{\frac{\partial }{\partial x}}+
\eta_{2}{\frac{\partial }{\partial y}}+
\eta_{2}^{(1)}{\frac{\partial }{\partial y'}}\:,
\end{eqnarray}
and whose commutator is given by\footnote{For the cases where $\left[X_{1},X_{2}\right] \neq 0$ see, for example, \cite{BluAnc}.} 
\begin{equation}
\left[X_{1},X_{2}\right]=0.
\end{equation}

\item Find the first integrals: The symmetries $\:X_{1}\:$ and $\:X_{2}\:$ are such that
\begin{equation}
\Delta = \left|
\begin{array}{ccc}
1 & y' & \phi \\
\xi_{1} & \eta_{1} & \eta_{1}^{(1)} \\
\xi_{2} & \eta_{2} & \eta_{2}^{(1)}
\end{array}
\right| \neq 0,
\end{equation}
and so, we always have two functions $\:I_1(x,y,y')\:$ and $\:I_2(x,y,y')\:$
satisfying 
\begin{eqnarray}
D_x[{I_1}] & = & {I_1}_{x} + y'{I_1}_{y} + \phi{I_1}_{y'}=0, \label{sistphi1}\\
X_{1}[{I_1}]&=&\xi_{1}{I_1}_{x}+\eta_{1}{I_1}_{y}
+\eta_{1}^{(1)}{I_1}_{y'}=0, \label{sistphi2}\\
X_{2}[{I_1}]&=&\xi_{2}{I_1}_{x}+\eta_{2}{I_1}_{y}
+\eta_{2}^{(1)}{I_1}_{y'}=1,
\label{sistphi3}
\end{eqnarray}
and
\begin{eqnarray}
D_x[{I_2}] & = & {I_2}_{x} + y'{I_2}_{y} + \phi{I_2}_{y'}=0, \label{sistpsi1}\\
X_{1}[{I_2}]&=&\xi_{1}{I_2}_{x}+\eta_{1}{I_2}_{y}+\eta_{1}^{(1)}{I_2}_{y'}=1, \label{sistpsi2}\\
X_{2}[{I_2}]&=&\xi_{2}{I_2}_{x}+\eta_{2}{I_2}_{y}+\eta_{2}^{(1)}{I_2}_{y'}=0. 
\label{sistpsi3}
\end{eqnarray}
Since the determinants of (\ref{sistphi1},\ref{sistphi2},\ref{sistphi3}) and (\ref{sistpsi1},\ref{sistpsi2},\ref{sistpsi3}) are $\neq 0$, the systems can be solved for
$\:{I_1}_{x}, {I_1}_{y}, {I_1}_{y'}\:$ and 
$\:{I_2}_{x}, {I_2}_{y}, {I_2}_{y'}\:$. So $\:{I_1}(x,y,y')\:$ and $\:{I_2}(x,y,y')\:$ can be written as line integrals.

\end{enumerate}

\begin{obs}
When the 2ODE does not present Lie point symmetries, there is no general method to look for the symmetry groups.
\end{obs}

\begin{obs}
We can deal with dynamical symmetries if the 2ODE presents two of them in the same way as above, but we can not apply the strategies above for the case of a single non point symmetry.
\end{obs}

\medskip

\subsection{The case of one non point Lie symmetry}
\label{npsym}

Supose that the 2ODE (\ref{r2ode}) admits a symmetry group given by
\begin{equation}
X^{(1)} = \xi(x,y,y')\,\frac{\partial}{\partial x}+\eta(x,y,y')\,\frac{\partial}{\partial y}+\zeta(x,y,y')\,\frac{\partial}{\partial y'}, \nonumber 
\end{equation}
where\footnote{Otherwise the transformations would not form a group.} $\zeta(x,y,y')=D_x\left[\eta(x,y,y')\right] - y'\,D_x\left[\xi(x,y,y')\right]$. In this case the group does not present an invariant of zeroth order. So, to use the differential invariants, we have to solve the characteristic system of the PDE
\begin{equation}
\xi(x,y,y')\,\frac{\partial \psi}{\partial x}+\eta(x,y,y')\,\frac{\partial \psi}{\partial y}+\zeta(x,y,y')\,\frac{\partial \psi}{\partial y'}=0
\end{equation}
that is given by
\begin{equation}
\frac{dx}{\xi(x,y,y')}=\frac{dy}{\eta(x,y,y')}=\frac{dy'}{\zeta(x,y,y')},
\end{equation}
which is, probably, more difficult to solve than the original 2ODE. In this case, what is usually done is to use the result described below to put the symmetry in evolutionary form:

\begin{teor} 
\label{ressydy}
If  $X^{(1)}=\xi \partial_{x}+\eta \partial_{y}+\eta^{(1)}
\partial_{y'}$ is a symmetry admitted by the 2ODE (\ref{r2ode}), then
\begin{equation}
\tilde{X}^{(1)} \equiv X^{(1)} + \rho(x,y,y') \, D_x
\end{equation}
is also a symmetry.
\end{teor}

\begin{prov}
\begin{eqnarray}
&&\left[D_x,\tilde{X}^{(1)}\right]=\left[D_x,X^{(1)} +\rho \,D_x\right]=\left[D_x,X^{(1)}\right]+\left[D_x,\rho D_x\right]= \nonumber \\ [2mm]
&&=D_x[\xi]\,D_x+D_x[\rho]\,\,D_x= D_x[\tilde{\xi}]\,D_x. \,\,\,\,\, \Box \nonumber
\end{eqnarray}
\end{prov}

\noindent 
From theorem \ref{ressydy}, by choosing $\:\rho=-\xi$, we have that $\:\tilde{X}\:$ will not have the coefficient of\footnote{$\partial_u \equiv \displaystyle{\frac{\partial}{\partial u}}$.} $\partial_{x}$:
\begin{equation}
\tilde{X}^{(1)}=\left(\eta-y' \xi\right)\partial_{y}+
\left(\eta^{(1)}-\phi \xi\right)\partial_{y'},
\end{equation}
and, besides this, the commutator $\:[D_x,\tilde{X}^{(1)}]\:=0$. So, with a symmetry in the {\em evolutionary form}\footnote{A symmetry with $\:\xi=0\:$ is said to be in {\em evolutionary form}.} we have that $I^{(0)}=x$. In this way, we can use the method of differential invariants.

\section{A method to search for symmetries}
\label{themet}

In this section we construct a method to find the symmetries of a rational 2ODE that presents (at least) one Liouvillian first integral. We present a formal way to represent a symmetry that helped us devise a method to look for it. Then we show that, although not entirely general, this algorithm can be applied to a large class of 2ODEs very difficult to treat (not only) by symmetry methods. After that, we present the steps of our algorithm and an example to clarify the `inner works' of our procedure.

\subsection{An useful way to represent a symmetry}
\label{formalsym}

Let us return to the rational 2ODE (\ref{r2ode}). Consider that ${X_e}^{(1)}$ defines a symmetry for (\ref{r2ode}) in the evolutionary form
\begin{equation}
\label{x1evol}
{X_e}^{(1)}=\nu\,\partial_y+D_x[\nu]\,\partial_{y'},
\end{equation}
\noindent
where $\,\nu \equiv \eta-y'\,\xi\,$. In this case the symmetry condition is ${X_e}^{(2)}\,[y''-\phi(x,y,y')] = 0, \, \mod(y''-\phi(x,y,y')=0)$, where ${X_e}^{(2)} = \nu\,\partial_y+D_x[\nu]\,\partial_{y'}+D_x^2[\nu]\,\partial_{y''}$. So, $\nu$ obeys the equation 
\begin{equation}
\label{eqsymev}
D_x^2[\nu] - D_x[\nu]\,\phi_{y'}-\nu\,\phi_y = 0,
\end{equation}
\noindent
which is a linear homogeneous second order PDE (linear homogeneous 2PDE). The fundamental idea (that culminated in this work) came from the formal equivalence between the operators $D_x$ and $\displaystyle{\frac{d}{dx}}$ over the solutions of the rational 2ODE (\ref{r2ode}). Because of this equivalence, over the solutions, the 2PDE (\ref{eqsymev}) can be regarded as a linear homogeneous 2ODE. Now, it is a well-known fact that any Riccati 1ODE is related (by a change of variables) with a linear homogeneous 2ODE. Let's see how it happens:
\begin{defin}[Riccati equation]
\label{riccati}
Any 1ODE that can be written as
\begin{equation}
y'=f(x)y^2+g(x)y+h(x),
\label{eqriccati}
\end{equation}
where $f(x)$, $g(x)$ and $h(x)$ are infinitely differentiable functions ($C^{\infty}$), is denominated {\bf Riccati equation} or {\bf Riccati 1ODE}.\footnote{These equations are the most general 1EDOs whose movable singularities are poles.}
\end{defin}
Making the following change of variables
\begin{equation}
y=-{\frac{1}{f(x)}}{\frac{w'}{w}}
\label{changericcati}
\end{equation}
in the Riccati equation (\ref{changericcati}), we obtain
\begin{equation}
w''={\frac {f'(x)+g(x)f(x)}{f(x)}}w'-h(x)f(x)w,
\label{eqriccatiL2}
\end{equation}
which is a linear homogeneous 2ODE. Thus, considering the 2PDE (\ref{eqsymev}) as a linear homogeneous 2ODE (over the solutions of the 2ODE (\ref{r2ode})), the analogous of the transformation (\ref{changericcati}) is 
\begin{equation}
\sigma = - \frac{D_x\left[\nu\right]}{\nu},
\label{changericcati2}
\end{equation}
and its inverse is given by the formal solution of (\ref{changericcati2}) for $\nu$:
\begin{equation}
\nu = {\rm e}^{-\int_x\left[\sigma\right]}.
\label{invchangericcati2}
\end{equation}

\begin{obs}
Note that the operator $D_x$ is a total derivative only over the solutions of the 2ODE. So, $\int_x$ is the inverse operator$\,$\footnote{In other words, $\int_x$ is a non local operator.} of $D_x$, i.e.,
\begin{equation}
\int_x\,D_x = D_x\,\int_x = \mbox{\boldmath $1$}.
\label{Dintx}
\end{equation}
\end{obs}

From the results above we can enunciate the following theorem:

\begin{teor}
\label{teosig}
If the rational 2ODE {\em (\ref{r2ode})} admits a symmetry given by ${X_e}^{(1)}=\nu\,\partial_y+D_x[\nu]\,\partial_{y'}$, then the function $\sigma$ defined by 
\begin{equation}
\sigma \equiv - \frac{D_x\left[\nu\right]}{\nu},
\label{defsig}
\end{equation}
obeys the 1PDE 
\begin{equation}
D_x\left[\sigma \right] = \sigma^2 +\partial_{y'}[\phi]\,\sigma - \partial_{y}[\phi].
\label{eqsig}
\end{equation}
Conversely, if the function $\,\sigma$ is a solution of the 1PDE (\ref{eqsig}) then the function $\,\nu \equiv {\rm e}^{-\int_x\left[\sigma\right]}$ defines a symmetry for the rational 2ODE {\em (\ref{r2ode})}.
\end{teor}

\begin{prov}
If the hypothesis of the theorem is satisfied, there exists a function $\nu$ satisfying $D_x^2[\nu] = D_x[\nu]\,\partial_{y'}[\phi] + \nu\,\partial_{y}[\phi]$. Substituting 
$\displaystyle{\nu = {\rm e}^{-\int_x\left[\sigma\right]}}$ in this 2PDE, since
\begin{eqnarray}
D_x\left[\nu\right] &=& D_x\left[{\rm e}^{-\int_x\left[\sigma\right]}\right] = {\rm e}^{-\int_x\left[\sigma\right]}\,D_x\left[-\int_x\left[\sigma\right]\right] = 
- \sigma\,{\rm e}^{-\int_x\left[\sigma\right]}, \nonumber \\ [2mm]
{D_x}^2\left[\nu\right] &=& D_x\left[- \sigma\,{\rm e}^{-\int_x\left[\sigma\right]}\right] = D_x\left[- \sigma\right]\,{\rm e}^{-\int_x\left[\sigma\right]} - \sigma\,D_x\left[{\rm e}^{-\int_x\left[\sigma\right]}\right]  =  \nonumber \\ [2mm]
&=& - D_x\left[\sigma\right]\,{\rm e}^{-\int_x\left[\sigma\right]} - \sigma\,\left(- \sigma\,{\rm e}^{-\int_x\left[\sigma\right]}\right),
\nonumber
\end{eqnarray}
we obtain
\begin{equation}
- D_x\left[\sigma\right]\,{\rm e}^{-\int_x\left[\sigma\right]} - \sigma\,\left(- \sigma\,{\rm e}^{-\int_x\left[\sigma\right]}\right) = 
- \sigma\,{\rm e}^{-\int_x\left[\sigma\right]}\,\partial_{y'}[\phi] + {\rm e}^{-\int_x\left[\sigma\right]}\,\partial_{y}[\phi] , \nonumber
\end{equation}
\begin{equation}
\Rightarrow \,\,\, D_x\left[\sigma \right] = \sigma^2 +\partial_{y'}[\phi]\,\sigma - \partial_{y}[\phi].  \nonumber
\end{equation}
Now, let's prove the converse: Consider the function defined by $\,\nu = {\rm e}^{-\int_x\left[\sigma\right]}$ where $\sigma$ is a solution of the 1PDE (\ref{eqsig}). So, 
\begin{eqnarray}
&& D_x^2[\nu] - D_x[\nu]\,\partial_{y'}[\phi] - \nu\,\partial_{y}[\phi] = \nonumber \\ [2mm]
&=& - D_x\left[\sigma\right]\,{\rm e}^{-\int_x\left[\sigma\right]} - \sigma\,\left(- \sigma\,{\rm e}^{-\int_x\left[\sigma\right]}\right) 
+\phi_{y'}\,\sigma\,{\rm e}^{-\int_x\left[\sigma\right]}-\phi_y\,\,{\rm e}^{-\int_x\left[\sigma\right]} = \nonumber \\ [2mm]
&=& {\rm e}^{-\int_x\left[\sigma\right]} \left( - D_x\left[\sigma\right] + \sigma^2+\phi_{y'}\,\sigma-\phi_y \right). \nonumber 
\end{eqnarray}
Since, by hypothesis, $\sigma$ is a solution of the 1PDE (\ref{eqsig}) then the term $- D_x\left[\sigma\right] + \sigma^2+\phi_{y'}\,\sigma-\phi_y$ equals zero. So, $D_x^2[\nu] - D_x[\nu]\,\partial_{y'}[\phi] - \nu\,\partial_{y}[\phi] =0. \,\,\,\,\Box$ 
\end{prov}

\medskip

\subsection{The relation with an integrating factor}
\label{relif}

Let us assume now that the 2ODE has a first integral that can be expressed in `closed form', that is, it presents a Liouvillian first integral. The knowledge of an integrating factor allows to obtain this first integral through quadratures. In the following we will show how this integrating factor is related to a symmetry (in evolutionary form).

\begin{defin}[Integrating factor]
\label{defif}
Let $\,y''=\phi(x,y,y')$ be a rational 2ODE presenting a Liouvillian first integral $I(x,y,y')$. Any function $\,\mu(x,y,y')$ satisfying
\begin{equation}
\mu(x,y,y')\,\left(\phi(x,y,y')-y''\right)=\frac{d}{dx} \left[I(x,y,y')\right],
\label{eqdefif}
\end{equation}
is denominated an {\bf integrating factor} for the 2ODE.
\end{defin}

\noindent
From (\ref{eqdefif}) we can write\footnote{$d_x \equiv \displaystyle{\frac{d}{dx}}$.}
\begin{equation}
\mu\,\phi-\mu\,y''= d_x \left[I(x,y,y')\right] = \partial_x[I]+y'\,\partial_y[I]+y''\,\partial_{y'}[I],
\label{relmui}
\end{equation}
So, since in (\ref{relmui}) neither $\mu$ nor $\phi$ nor $I$ are functions of $y''$, we have that
\begin{equation}
\mu=-\partial_{y'}[I].
\label{relmuI}
\end{equation}
On the other hand, suppose that ${X_e}^{(1)}=\nu\,\partial_y+D_x[\nu]\,\partial_{y'}$ is a symmetry vector field of the rational 2ODE $\,y''=\phi(x,y,y')$ such that
\begin{equation}
{X_e}^{(1)}[I]=\nu\,\partial_y[I]+D_x[\nu]\,\partial_{y'}[I]=0.
\label{xei}
\end{equation}
We can write (\ref{xei}) as
\begin{equation}
-\frac{D_x[\nu]}{\nu} =\frac{\partial_y[I]}{\partial_{y'}[I]},
\label{defsigm}
\end{equation}
and noting that $\displaystyle{-\frac{D_x[\nu]}{\nu}}$ is $\sigma$ (as defined in (\ref{defsig})), we have that
\begin{equation}
\sigma =\frac{\partial_y[I]}{\partial_{y'}[I]},
\label{relsigI}
\end{equation}

\noindent
Based on (\ref{relmuI}) and (\ref{relsigI}) we can enunciate the following result:

\begin{teor}
\label{relmusig}
Let $\,y''=\phi(x,y,y')$ be a rational 2ODE presenting a Liouvillian first integral $I(x,y,y')$. Also, let $\mu$ be an integrating factor such that $\mu=-\partial_{y'}[I]$ and let $\nu$ be a function defining a symmetry (in evolutionary form) such that ${X_e}^{(1)}[I]=\nu\,\partial_y[I]+D_x[\nu]\,\partial_{y'}[I]=0$. Then the following equation holds:
\begin{equation}
-\frac{D_x[\mu]}{\mu} = -\frac{D_x[\nu]}{\nu} + \partial_{y'} [\phi].
\label{relmusym}
\end{equation}
\end{teor}

\begin{prov}
\

\noindent
If the hypotheses of the theorem are satisfied, then $\mu=-\partial_{y'}[I]$. Applying $D_x$, we obtain
\begin{equation}
D_x[\mu] = -D_x\left[\partial_{y'} [I]\right].
\label{eqproo1}
\end{equation}
Since the commutator $\displaystyle{\left[\partial_{y'},D_x\right]}$ is given by
\begin{equation}
\left[\partial_{y'},D_x\right] = \partial_{y'}\,D_x - D_x\,\partial_{y'} = \partial_y + \partial_{y'}[\phi]\,\partial_{y'},
\label{eqproo2}
\end{equation}
then
\begin{eqnarray}
&& \left[\partial_{y'},D_x\right][I] = \partial_{y'}\,\left[D_x[I]\right] - D_x\,\left[\partial_{y'}[I]\right] = \partial_y[I] + \partial_{y'}[\phi]\,\partial_{y'}[I] = \nonumber \\ [2mm]
&& \partial_{y'}\big[\,\overbrace{D_x[I]}^{=\:0}\big] + \overbrace{D_x\,\left[-\partial_{y'}[I]\right]}^{=\:D_x[\mu]} = \partial_y[I] - \partial_{y'}[\phi]\,\mu.
\label{eqproo2}
\end{eqnarray}
Dividing (\ref{eqproo2}) by $\partial_{y'}[I]$ and observing equation (\ref{defsigm}), we obtain finally
\begin{eqnarray}
&& \frac{D_x[\mu]}{\partial_{y'}[I]} = \frac{\partial_y[I]}{\partial_{y'}[I]} - \partial_{y'}[\phi]\,\frac{\mu}{\partial_{y'}[I]} 
\nonumber \\ [2mm]
&\Rightarrow& -\frac{D_x[\mu]}{\mu} = -\frac{D_x[\nu]}{\nu} + \partial_{y'}[\phi]. \,\,\,\, \Box
\label{eqproo3}
\end{eqnarray}
\end{prov}

\subsection{A possible algorithm}
\label{possalg}

In this paper, we are interested in finding symmetries to use them in the same old spirit of Lie's theory, that is, to integrate 2ODEs using quadratures. In this sense, when talking about rational 2ODEs with a Liouvillian first integral, we are looking for integrating factors of the form\footnote{We could not actually find any case where this was not valid. So, we would say that the vast majority of cases are covered by this structure.}
\begin{equation}
\mu = {\rm e}^{Z_0}\,\prod_i {p_i}^{n_i},
\label{strucif}
\end{equation}
where $Z_0$ is a rational function of $(x,y,y')$, the $p_i$ are irreducible polynomials in $(x,y,y')$ and the $n_i$ are rational numbers. In this case, we can formulate the following theorem:

\begin{teor}
\label{relmusigrat}
Let $\,y''=\phi(x,y,y')$ be a rational 2ODE presenting a Liouvillian first integral $I(x,y,y')$ and let $\mu$ be an integrating factor of the form  
\begin{equation}
\mu = {\rm e}^{Z_0}\,\prod_i {p_i}^{n_i},
\label{strucif}
\end{equation}
where $Z_0$ is a rational function of $(x,y,y')$, the $p_i$ are irreducible polynomials in $(x,y,y')$, the $n_i$ are constants, such that $\mu=-\partial_{y'}[I]$. Also, let $\nu$ be a function defining a symmetry (in evolutionary form) such that ${X_e}^{(1)}[I]=\nu\,\partial_y[I]+D_x[\nu]\,\partial_{y'}[I]=0$. Then the function $\sigma$ defined by
\begin{equation}
\sigma \equiv - \frac{D_x\left[\nu\right]}{\nu},
\label{defsignu}
\end{equation}
(where $D_x$ is the vector field associated with the 2ODE), is a rational function of $(x,y,y')$.
\end{teor}

\begin{prov}
\

\noindent
Assuming that $\mu = {\rm e}^{Z_0}\,\prod_i {p_i}^{n_i}$, we have 
\begin{equation}
\frac{D_x[\mu]}{\mu} = \frac{ e^{Z_0}\,D_x\left[ \prod_i p_i^{n_i} \right] + e^{Z_0}\,D_x\left[ Z_0 \right] \, \prod_i p_i^{n_i} }{e^{Z_0}\, \prod_i p_i^{n_i}} = D_x\left[ Z_0 \right] + \frac{D_x\left[ \prod_i p_i^{n_i} \right] }{ \prod_i p_i^{n_i}}. \nonumber
\end{equation}
Since $\,D_x = \partial_x + y'\,\partial_y + \phi(x,y,y')\,\partial_{y'}\,$, $\,Z_0\,$ is a rational function of $(x,y,y')$ and the $p_i$ are polynomials in $(x,y,y')$, then we have that $\frac{D_x[\mu]}{\mu}$ is a rational function of $(x,y,y')$. However, from the result of theorem \ref{relmusig}
\begin{equation}
-\frac{D_x[\mu]}{\mu} = -\frac{D_x[\nu]}{\nu} + \partial_{y'} [\phi] = \sigma + \partial_{y'} [\phi].
\label{relmusym}
\end{equation}
Finally, since $\phi$ is rational, so is $\partial_{y'} [\phi]$. Therefore, $\sigma$ is a rational function of $(x,y,y')$.  $\,\,\Box$
\end{prov}

\begin{obs}
Note that, although $\mu$ is an elementary function of the form shown above, the symmetry $\nu$ need not be of such format, i.e., all we have is that $\,\nu = {\rm e}^{-\int_x\left[\sigma\right]}$, where $\sigma$ is a rational function. On the other hand, if $\nu$ were an elementary function of the format shown above, then $\sigma$ would also be a rational function, but now the integrating factor would not necessarily be an elementary function.
\end{obs}

From theorems \ref{teosig} and \ref{relmusigrat}, we can establish a result that allows for the construction of an algorithm to search for symmetries of rational 2ODEs:

\begin{teor}
\label{relmusigrat}
Let $$\,y''=\phi(x,y,y')=\displaystyle{\frac{M(x,y,y')}{N(x,y,y')}},$$ where $M$ and $N$ are co-prime polynomials, be a rational 2ODE presenting a global analytical first integral $I(x,y,y')$ and a symmetry vector field ${X_e}^{(1)}=\nu\,\partial_y+D_x[\nu]\,\partial_{y'}$ such that the function $\sigma$ defined as $$\sigma \equiv \displaystyle{-\frac{D_x[\nu]}{\nu}}$$ is a rational function, i.e., $\sigma={\frac{p}{q}}$, where $p$ and $q$ are co-prime polynomials. Let $deg_M$, $deg_N$, $deg_p$ and $deg_q$ be the degrees of polynomials $M$, $N$, $p$ and $q$, respectively. If $deg_M \leq deg_N + 1$ then $deg_p \leq deg_q$. On the contrary, if $deg_M > deg_N + 1$ then $deg_p \leq deg_q + deg_M - deg_N - 1.$
\end{teor}

\begin{prov}
\

\noindent
From theorem \ref{teosig}, the function $\sigma$ obeys the equation 
$$
D_x\left[\sigma \right] = \sigma^2 +\partial_{y'}[\phi]\,\sigma - \partial_{y}[\phi].
$$
Substituting $\sigma={\frac{p}{q}}$, we have
$$
\frac{q\,D_x\left[p \right]-p\,D_x\left[q \right]}{q^2} = \frac{p^2}{q^2} +\partial_{y'}[\phi]\,\frac{p}{q} - \partial_{y}[\phi],
$$
implying that\footnote{Where $v_u$ means $\partial_u [v]$.}
\begin{equation}
\label{eqproofpq1}
\frac{p_x}{q}+\frac{y'\,p_y}{q}+\frac{M\,p_{y'}}{N\,q} - \frac{p\,q_x}{q^2} - \frac{y'\,p\,q_y}{q^2} - \frac{M\,p\,q_{y'}}{N\,q^2}
= \frac{p^2}{q^2} + \frac{p\,M_{y'}}{q\,N} - \frac{p\,M\,N_{y'}}{q\,N^2} - \frac{M_{y}}{N} + \frac{M\,N_{y}}{N^2}.
\end{equation}
Multiplying the equation (\ref{eqproofpq1}) by $q^2N^2$ and isolating the term $p^2N^2$ we obtain
$$
p^2N^2 = - p\,q\,N\,M_{y'} + p\,q\,M\,N_{y'} + q^2N\,M_{y} - q^2M\,N_{y} + 
$$
\begin{equation}
+ q\,p_xN^2 + y'q\,p_yN^2 + q\,p_{y'}N\,M - p\,q_xN^2 - y'p\,q_yN^2 - p\,q_{y'}N\,M. \label{eqproofpq2}
\end{equation}
The degree of the term in the left hand side of equation (\ref{eqproofpq2}) is $\,2\,deg_p + 2\,deg_N,\,$ since it is a square ($p^2N^2$). The (maximum) degrees of the following terms (in the right hand side) are, respectively, 
\begin{eqnarray}
&& deg_p + deg_q + deg_M + deg_N - 1, \nonumber \\ 
&& deg_p + deg_q + deg_M + deg_N - 1, \nonumber \\ 
&& 2\,deg_q + deg_M + deg_N - 1, \nonumber \\ 
&& 2\,deg_q + deg_M + deg_N - 1, \nonumber \\ 
&& deg_p + deg_q + 2\,deg_N - 1, \nonumber \\ 
&& deg_p + deg_q + 2\,deg_N, \nonumber \\ 
&& deg_p + deg_q + deg_M + deg_N - 1, \nonumber \\ 
&& deg_p + deg_q + 2\,deg_N - 1, \nonumber \\ 
&& deg_p + deg_q + 2\,deg_N, \nonumber \\ 
&& deg_p + deg_q + deg_M + deg_N - 1. \nonumber
\end{eqnarray}
There are four cases:

\noindent
{\bf case 1:} $2\,deg_p + 2\,deg_N \leq deg_p + deg_q + deg_M + deg_N - 1$;

\noindent
{\bf case 2:} $2\,deg_p + 2\,deg_N \leq 2\,deg_q + deg_M + deg_N - 1$;

\noindent
{\bf case 3:} $2\,deg_p + 2\,deg_N \leq deg_p + deg_q + 2\,deg_N - 1$;

\noindent
{\bf case 4:} $2\,deg_p + 2\,deg_N \leq deg_p + deg_q + 2\,deg_N$;

\vspace{2mm}
\noindent
So, one of the following inequalities holds:

\begin{enumerate}
\item $deg_p \leq deg_q + deg_M - deg_N - 1$;
\item $deg_p \leq deg_q + \frac{deg_M - deg_N - 1}{2}$;
\item $deg_p \leq deg_q - 1$;
\item $deg_p \leq deg_q$;
\end{enumerate}

\noindent
The third inequality implies the fourth inequality and second inequality implies the first one. So, if $deg_M > deg_N + 1$ then $deg_p \leq deg_q + deg_M - deg_N - 1$, otherwise, if $deg_M \leq deg_N + 1$ then $deg_p \leq deg_q$. $\,\,\,\Box$

\end{prov}

\vspace{2mm}
From the result above, we can establish a semi algorithm\footnote{An algorithm (in the sense we are considering here) is a set of rules that can be translated into a programming language so that the rule set can be executed by a computer in a finite number of steps. Our procedure (in the general case) is a semi-algorithm because it is not possible to determine a priory the degree of the rational function $\sigma$.} to find the symmetries of a rational 2ODE. The idea is to use the equation (\ref{eqproofpq2}) as the basis of our procedure: we basically construct two polynomials with undetermined coefficients (candidates for $p$ and $q$) and substitute them in equation (\ref{eqproofpq2}). Then we collect the resulting polynomial equation in monomials of $(x,y,y')$ and equal each of the coefficients to zero. By doing so, we will find a system of algebraic equations for the coefficients that, if solved, will allow us to formally write the following symmetry:
\begin{equation}
{X_e}^{(1)} =\nu\,\partial_y+D_x[\nu]\,\partial_{y'}={\rm e}^{-\int_x\left[\sigma\right]}\,\left(\partial_y - \sigma\,\partial_{y'}\right),
\label{xei2}
\end{equation}
where $\,\sigma\,$ is constructed with the solution ($p/q$) of the system of algebraic equations. Although the symmetry is written in terms of the functional $\displaystyle{\int_x\left[\sigma\right]}$, its use for the calculation of the invariants $I^{(0)}$ and $I^{(1)}$ does not depend on it: $I^{(0)}=x$ (evolutionary form) and $I^{(1)}$ is such that $\displaystyle{{\rm e}^{-\int_x\left[\sigma\right]}\,\left(\partial_y - \sigma\,\partial_{y'}\right)\left[I^{(1)}\right]=0}$, implying that  $\displaystyle{\left(\partial_y - \sigma\,\partial_{y'}\right)\left[I^{(1)}\right]=0}$. So, we can use the method of differential invariants in the same way we did in section \ref{npsym}.

\medskip

\subsubsection{The steps of the (semi) algorithm}
\label{stepsa}

The last paragraphs of the previous section describe (colloquially) a possible (semi) algorithm for finding the symmetries of a rational 2ODE.

\begin{algor}  (ASymm)
\label{asymm}
	\begin{enumerate}
           \item Construct the $D_x$ operator.
		\item Set $n_{max}$={\tt some integer}, $d_M=degree(M)$ and  $d_N=degree(N)$.
		\item Set $n=0$.
		\item Set $n=n+1$.
		\item if $n>n_{max}$ then FAIL.
           \item {\bf If}  $d_M > d_N + 1$ {\bf then} construct a generic polynomial ${q_c}$ of degree $n$ in $(x,y,y')$ with undetermined coefficients $b_i$ and a generic polynomial ${p_c}$ of degree $n + d_M - d_N - 1$ in $(x,y,y')$ with undetermined coefficients $a_i$ {\bf else}  construct two generic polynomials ${p_c}$ and ${q_c}$ of degree $n$ in $(x,y,y')$ with undetermined coefficients $a_i$ and $b_i$, respectively.
		\item Substitute $p_c$ and ${q_c}$ in equation (\ref{eqproofpq2}), collect the resulting polynomial equation in the variables $x,\,y,\,y'$ and equate the coefficients of each monomial to zero, obtaining a system $AE$ of algebraic equations.
		\item Solve the system ${AE}$ with respect to $a_i$ and $b_i$. If no solution is found, then go to step 4.
		\item Substitute the solution in $p_c/q_c$ (obtaining $\sigma$).
	\end{enumerate}
\end{algor}

\subsubsection{Example}
\label{examp}

Consider the following rational 2ODE:
\begin{equation}
\label{ex1-2ode}
y''={\frac { \left( y'-1 \right)  \left( {x}^{4}y'+2\,{x}^{3}y-{x}^{2}y+y'
 \right) }{ \left( {x}^{2}y-1 \right) {x}^{2}}}.
\end{equation}
Let's apply the procedure steps to the 2ODE (\ref{ex1-2ode}):\footnote{In this example we show in detail the steps of the algorithm $ASymm$.}
\begin{enumerate}
\item  $D_x = \partial_x+y'\partial_y+{\frac { \left( y'-1 \right)  \left( {x}^{4}y'+2\,{x}^{3}y-{x}^{2}y+y'
 \right) }{ \left( {x}^{2}y-1 \right) {x}^{2}}}\,\partial_{y'}.$
\item $n_{max}=7$, $d_M=6$ and  $d_N=5$.
\item $n=0$.
\item $n=n+1=1$.
\item $1<7$ (proceed).
\item $d_M = d_N +1 \,\,\, \Rightarrow$ 

$p_c=a_{1}\,x+a_{2}\,y+a_{3}\,y'+a_{0}$, 

$q_c=b_{1}\,x+b_{2}\,y+b_{3}\,y'+b_{0}$.

\item The system $AE$ has only the trivial solution (go to step 4).

\item $n=n+1=2$.
\item $2<7$ (proceed).
\item 

$p_c=a_{4}\,{x}^{2}+a_{5}\,{y}^{2}+a_{6}\,{y'}^{2}+a_{7}
\,xy'+a_{8}\,yx+a_{9}\,y'y+a_{1}\,x+a_{2}\,y+a_{3}\,y'+a_{0}$, 

$q_c=b_{4}\,{x}^{2}+b_{5}\,{y}^{2}+b_{6}\,{y'}^{2}+b_{7}
\,xy'+b_{8}\,yx+b_{9}\,y'y+b_{1}\,x+b_{2}\,y+b_{3}\,y'+b_{0}$.

\item The system $AE$ has only the trivial solution (go to step 4).

\item $n=n+1=3$.
\item $3<7$ (proceed).
\item 

$p_c=a_{0}+a_{1}\,x+a_{2}\,y+a_{3}\,y'+a_{4}\,{x}^{
2}+a_{5}\,{x}^{3}+a_{6}\,{y}^{2}+a_{7}\,{y}^{3}+a_{8}\,{y'}^{2}+a_{9}\,{y'}^{3}+a_{10}\,x{y}^{2}+a_{11}\,
xy'+a_{12}\,x{y'}^{2}+a_{13}\,{x}^{2}y+a_{14}\,{x}^{2}y'+a_{15}\,yx+a_{16}\,y{y'}^{2}+a_{17}\,{y}^{2}y'+a_{18}
\,y'y+a_{19}\,xyy'$, 

$q_c=b_{0}+b_{1}\,x+b_{2}\,y+b_{3}\,y'+b_{4}\,{x}^{
2}+b_{5}\,{x}^{3}+b_{6}\,{y}^{2}+b_{7}\,{y}^{3}+b_{8}\,{y'}^{2}+b_{9}\,{y'}^{3}+b_{10}\,x{y}^{2}+b_{11}\,
xy'+b_{12}\,x{y'}^{2}+b_{13}\,{x}^{2}y+b_{14}\,{x}^{2}y'+b_{15}\,yx+b_{16}\,y{y'}^{2}+b_{17}\,{y}^{2}y'+b_{18}
\,y'y+b_{19}\,xyy'$.

\item The system $AE$ has the solution 

$a_{0}=0,a_{1}=0,a_{10}=0,a_{11}=0,a_{12}=0,a_{13}=0,a_{14}=b_{0},a_{15}=0,a_{16
}=0,a_{17}=0,a_{18}=0,a_{19}=0,a_{2}=0,a_{3}=0
,a_{4}=-b_{0},a_{5}=0,a_{6}=0,a_{7}=0,a_{8}=0,a_{9}=0,b_{0}=b_{0},b_{1}=0,b_{10}=0,
b_{11}=0,b_{12}=0,b_{13}=-b_{0},b_{14}=0,b_{15}=0,b_{16}=0,b_{17}=0,b_{18}=0,b_{19}=0,
b_{2}=0,b_{3}=0,b_{4}=0,b_{5}=0,b_{6}=0,b_{7}=0,b_{8}=0,b_{9}=0.$

\item $\sigma=\displaystyle{-{\frac {{x}^{2} \left( y'-1 \right) }{{x}^{2}y-1}}}$.

\end{enumerate}

So, the 2ODE (\ref{ex1-2ode}) admits the symmetry 
\begin{equation}
{X_e}^{(1)} = {\rm e}^{-\int_x\left[-{\frac {{x}^{2} \left( y'-1 \right) }{{x}^{2}y-1}}\right]}\,\left(\partial_y + {\frac {{x}^{2} \left( y'-1 \right) }{{x}^{2}y-1}}\,\partial_{y'}\right).
\label{xeex1}
\end{equation}

\begin{obs}
Note that, although the 2ODE (\ref{ex1-2ode}) has no local symmetries, we can, using the method of differential invariants (see section \ref{npsym}), find a first integral: $${\frac {{{\rm e}^{\frac{1}{x}}} \left( {x}^{2}y-y' \right) }{y'-1}}$$
\end{obs}

\begin{obs}
The method we present is a semi algorithm because\footnote{A full algorithm must end in a finite number of steps.} we can not determine the maximum degree of the polynomials $p$ and $q$ that form the $\sigma$ function. However, up to the degree analyzed, we can be sure that if we do not find the $\sigma$ function, then it does not exist within the range considered. 
\end{obs}

\section{Interesting cases of the semi algorithm}
\label{speci}

In the last section we presented a method that can compute Lie symmetries of a rational 2ODE semi algorithmically. However, the method described above can be greatly improved in some special (but quite comprehensive) cases. 
These improvements can (depending on the 2ODE in question) greatly simplify the $ASymm$ semi algorithm. For instance, most Computer Algebra Systems (CAS) use Gr\"obner bases to solve systems of polynomial equations and, since the complexity of the problem grows exponentially with the increase in the number of coefficients, any decrease in the number of coefficients (and/or nonlinear terms) may represent the difference between success and failure. 

In this section we will discuss some of these possible improvements. 

\begin{obs}
At some points in this section and the next, times will be assigned to certain sets of operations. What does that mean? These times result from the simple execution of the steps (of these operations) in the CAS Maple 17 (In this paper all the computational data, time of running etc, was obtained on the same computer with the following configuration: Intel(R) Core(TM) i5-8265U @ 1.8 GHz).
\end{obs}

\subsection{Some general considerations}
\label{Sgc}
	 
In first place let's remember that we are looking for a rational function $\sigma=\frac{p}{q}$ that obeys the 1PDE $D_x\left[\sigma \right] = \sigma^2 +\partial_{y'}[\phi]\,\sigma - \partial_{y}[\phi]$, where $D_x \equiv \partial_{x} + y'\,\partial_{y} + \phi(x,y,y')\,\partial_{y'}$. So, if $\sigma$ is not a function of one of the variables $(x,y,y')$ then we can separate the 1PDE above into an overdetermined system of PDEs and, in this case, there is already an algorithm\footnote{See, for example, \cite{BluAnc}.} to determine its solution. Therefore, in what follows, we will consider that (at least in principle) $\sigma$ is a function of the three variables $(x,y,y')$. With this in mind, let's divide equation (\ref{eqproofpq2}) by $\,M\,N\,$, obtaining
$$
\frac{p^2N + p\,q\,M_{y'} - q^2M_{y} - q\,p_xN - y'q\,p_yN + p\,q_xN + y'p\,q_yN}{M} + 
$$
\begin{equation}
 + q\,\frac{q\,N_{y}- p\,N_{y'}}{N} - q\,p_{y'} + p\,q_{y'}=0. \label{eqpq4}
\end{equation}
A closer look at equation (\ref{eqpq4}) reveals a number of improvements that can be made to the $ASymm$ algorithm for certain particular situations: the term $\,- q\,p_{y'} + p\,q_{y'}\,$ is polynomial since $p$ and $q$ are polynomials. Besides that, this term can not be null because it is equal to $\,-q^2\partial_{y'}[\sigma]$ and, by assumption, neither $q$ nor $\partial_{y'}[\sigma]$ can be null. Since $M$ and $N$ are coprime, the terms 
\begin{equation}
\label{ovm}
\frac{p^2N + p\,q\,M_{y'} - q^2M_{y} - q\,p_xN - y'q\,p_yN + p\,q_xN + y'p\,q_yN}{M}
\end{equation}
and
\begin{equation}
\label{ovn}
q\,\frac{q\,N_{y}- p\,N_{y'}}{N} 
\end{equation}
are (each one of then) polynomials. In what follows, we will discuss some particular cases that allow us to improve the (semi) algorithm presented in section \ref{possalg}. 

\begin{obs}
As we commented earlier, the different improvements are directly linked with the format and degree of the candidates $p_c$ and $q_c$. Thus, to avoid very redundant presentations, we will only describe (in each of the following cases) the change in the format (and degree) of the $p_c$ and $q_c$ polynomials (we will represent the untouched steps of the algorithm with `` ... '' ).
\end{obs}

\medskip

\subsection{$N$ and $q$ do not present common factors}
\label{ndifeq}

To address this case, let's focus on the fact that the term 
\begin{equation}
\label{ovn}
q\,\frac{q\,N_{y}- p\,N_{y'}}{N} 
\end{equation}
(see equation (\ref{eqpq4})) is a polynomial. This leads to the following result:

\begin{teor}
\label{teonov}
Let $q,\,p$ and $N$ be polynomials defined as above. Consider that $\displaystyle{q\,\frac{q\,N_{y}- p\,N_{y'}}{N}}$ is a polynomial. If $N$ and $q$ do not present common polynomial factors then one the following conditions hold:
\begin{enumerate}
\item Neither $N_y$ nor $N_{y'}$ are null and $N\, | \, q\,N_{y}- p\,N_{y'}$.

\item $N$ is a function (only) of $x$.

\item $N_y=0$, $N_{y'} \neq 0$ and $N\, | \, p\,N_{y'}$. 
\end{enumerate}
\end{teor}

\begin{prov}
\

1. : Since $N$ and $q$ do not present common factors and neither $N_y$ nor $N_{y'}$ are null, then $\displaystyle{\frac{q\,N_{y}- p\,N_{y'}}{N}}$ has to be a polynomial, i.e., $N\, | \, q\,N_{y}- p\,N_{y'}$.

2. : $N$ is a function (only) of $x$ $\Rightarrow \,\, N_y=N_{y'}=0$. The whole term is null.

3. : Follows directly from substituting $N_y=0$ in $\displaystyle{q\,\frac{q\,N_{y}- p\,N_{y'}}{N}}$. The other apparent case $N_{y} \neq 0$, $N_{y'}=0$ is not possible because $N$ cannot divide $N_{y'}$. $\,\,\,\Box$

\end{prov}

\medskip

The first condition allows for using the fact that $q\,N_{y}- p\,N_{y'}=N\,P$ (where $P$ is a polynomial) for some `pre-calculations' of the coefficients of $p$ and $q$. These conditions can significantly reduce the difficulty of solving the algebraic equations for the coefficients of $p$ and $q$ that were extracted from the equation (\ref{eqproofpq2}) as a whole.

\medskip

The second condition ($N=N(x)$) implies that 
\begin{equation}
q\,\frac{p\,M_{y'} - q\,M_{y}}{M} + N\,\frac{p^2 - q\,p_x - y'q\,p_y + p\,q_x + y'p\,q_y}{M} - q\,p_{y'} + p\,q_{y'}=0. \label{eqnnx}
\end{equation}
Multiplying (\ref{eqnnx}) by $\frac{M}{N\,q}$ we obtain:
\begin{equation}
\frac{p\,M_{y'} - q\,M_{y}- M\,p_{y'}}{N} + p\,\frac{p + q_x + y'\,q_y}{q} - p_x - y'\,p_y - \frac{p\,q_{y'}M}{N\,q}=0. \label{eqnnx2}
\end{equation}
The term $- p_x - y'\,p_y$ is a polynomial or zero. In both cases we can improve our algorithm: for example, in the case where $- p_x - y'\,p_y$ is not zero, we have that the terms 
\begin{equation}
\frac{p\,M_{y'} - q\,M_{y}- M\,p_{y'}}{N}, \,\,\,\,\frac{p + q_x + y'\,q_y}{q}\,\,\,\,{\rm and} \,\,\,\,\frac{p\,q_{y'}M}{N\,q}\label{eqterms}
\end{equation}
must combine to form a polynomial and, since (by hypothesis) $N$ and $q$ do not present common factors we obtain that
\begin{equation}
\frac{p\,M_{y'} - q\,M_{y}- M\,p_{y'}-p\,q_{y'}M}{N}\,\,\,\,{\rm and} \,\,\,\,\frac{N\,(p + q_x + y'\,q_y)-p\,q_{y'}M}{q}
\label{eqterms2}
\end{equation}
are polynomials. A very likely case is that $q_{y '}$ is null and if so, this fact greatly simplifies the algorithm. 

\medskip

The third condition: $N_y=0$ and $N$ and $p$ have common factors, again would greatly simplify the algorithm because, since $N$ cannot divide $N_y$, the only way for $N_y$ and $N$ to have common factors is in the case that $N$ has polynomial factors raised to an integer exponent greater than 1. Otherwise, $N$ divides $p$.

\medskip

\subsubsection{$ASymm_{[1.1]}$}
\label{asymm11}

This subcase is linked to the condition: 

\medskip

\noindent
{\bf The term $\displaystyle{\frac{q\,N_{y}- p\,N_{y'}}{N}}$ is itself a polynomial (neither $N_y$ nor $N_{y'}$ are null)}.

\begin{algor}  ($ASymm_{[1.1]}$)
\label{asymm11alg}
	\begin{itemize}
           \item $\cdots$
           \item {\bf 6.}  {\bf If}  $d_M > d_N + 1$ {\bf then} construct a generic polynomial ${q_c}$ of degree $n$ in $(x,y,y')$ with undetermined coefficients $b_i$ and a generic polynomial ${p_c}$ of degree $n + d_M - d_N - 1$ in $(x,y,y')$ with undetermined coefficients $a_i$ {\bf else}  construct two generic polynomials ${p_c}$ and ${q_c}$ of degree $n$ in $(x,y,y')$ with undetermined coefficients $a_i$ and $b_j$, respectively.
           \item {\bf 6.1.}  Construct a generic polynomial ${P_c}$ of degree $\max(deg_{q_c},deg_{p_c})$ in $(x,y,y')$ with undetermined coefficients $c_k$. 
           \item {\bf 6.2.}  Substitute ${p_c}$, ${q_c}$ and ${P_c}$ in the equation $q_c\,N_{y}- p_c\,N_{y'}-N\,P_c=0$ and collect the resulting polynomial equation in the variables $x,\,y,\,y'$ and equate the coefficients of each monomial to zero, obtaining a system $AE_{[1.1]}$ of algebraic equations.
           \item {\bf 6.3.}  Solve the system $AE_{[1.1]}$ with respect to $a_i$ and $b_i$ (in terms of the $c_k$).
           \item {\bf 6.4.}  Substitute the result in the polynomials ${p_c}$ and ${q_c}$.
           \item $\cdots$
	\end{itemize}
\end{algor}

\begin{obs}
The solution obtained in the step 6.3 is the main point of this particular sub-algorithm ($ASymm_{[1.1}$). In practice, the number of undetermined coefficients   can greatly decrease helping the success of step 8.
\end{obs}

\subsubsection{$ASymm_{[1.2]}$}
\label{asymm12}

This subcase is linked to the condition: 

\medskip

\noindent
{\bf $N$ is a function (only) of $x$ (and $q_{y'}=0$)}.

\begin{algor}  ($ASymm_{[1.2]}$)
\label{asymm12alg}
	\begin{itemize}
           \item $\cdots$
           \item {\bf 6.}  {\bf If}  $d_M > d_N + 1$ {\bf then} construct a generic polynomial ${q_c}$ of degree $n$ in $(x,y)$ with undetermined coefficients $b_i$ and a generic polynomial ${p_c}$ of degree $n + d_M - d_N - 1$ in $(x,y,y')$ with undetermined coefficients $a_i$ {\bf else}  construct two generic polynomials ${p_c}$ and ${q_c}$ of degree $n$ in $(x,y,y')$ and  in $(x,y)$ with undetermined coefficients $a_i$ and $b_j$, respectively.
           \item $\cdots$
	\end{itemize}
\end{algor}

Let's take (for example) the 2ODE (nonlinear 2ODE 41 with $f(x)=1/x^2$ of Kamke's book \cite{Kam})
\begin{equation}
\label{2odeKam41}
y'' = -{\frac {{x}^{2}{y}^{3}-3\,{x}^{2}yy'+{y}^{2}-y'}{{x}^{2}}}
\end{equation}
with $N=N(x)$, has a $\sigma$ function given by
\begin{equation}
\label{sigKam41}
\sigma = -{\frac {2\,x{y}^{2}+2\,xy'+{y}^{2}+4\,y+y'}{2\,xy+y+2}},
\end{equation}
where the fact that $q_{y'}=0$ improves algorithm efficiency by more than 10 times.

\subsubsection{$ASymm_{[1.3]}$}
\label{asymm13}

This subcase is linked to the condition: 

\medskip

\noindent
{\bf $N_y=0$ and $N$ and $p$ have common factors}.

\begin{algor}  ($ASymm_{[1.3]}$)
\label{asymm13alg}
	\begin{itemize}
           \item $\cdots$
           \item {\bf 6.}  {\bf If}  $d_M > d_N + 1$ {\bf then} construct a generic polynomial ${q_c}$ of degree $n$ in $(x,y,y')$ with undetermined coefficients $b_i$ and a generic polynomial ${pp_c}$ of degree $n + d_M - d_N - 1 - d_f$, where $d_f$ is the degree of the common factors (CF) of $N$ and $p$ in $(x,y,y')$ with undetermined coefficients $a_i$; $p_c=pp_c\,CF$  {\bf else}  construct two generic polynomials ${pp_c}$ and ${q_c}$ of degree $n-d_f$ and $n$ in $(x,y,y')$, with undetermined coefficients $a_i$ and $b_j$, respectively; $p_c=pp_c\,CF$.
           \item $\cdots$
	\end{itemize}
\end{algor}

\subsection{$N$ and $q$ present common factors}
\label{neqq}

Regarding improving the efficiency of our semi algorithm it is more practical to divide the different situations into one or another of the two following distinct cases:

\begin{itemize}
\item $q$ has factors not present in $N$.
\item All factors of $q$ are factors of $N$.
\end{itemize}

\subsubsection{$q$ has factors not present in $N$}
\label{qdifern}

For this situation we can (for example) use $q_c = N q_m$ or $q_c = N_f\, q_m$ where $q_m$ has a lower degree and $N_f$ is a factor (or a product os some factors) of $N$. For example, the 2ODE (see section \ref{perfitK} 2ODE 180)
\begin{equation}
\label{2ode180}
y'' = {\frac {2\,({x}^{2}{y'}^{2}+yy'x-{y}^{3}-xy'+2\,{y}^{2}-y)}{{x}^{2} \left( y-1 \right) }}
\end{equation}
has a $\sigma$ given by 
\begin{equation}
\label{sigm180}
\sigma = -{\frac {y' \left( 2\,y-1 \right) }{y \left( y-1 \right) }}
\end{equation}
and, in order to determine it we do not need a second degree polynomial candidate, but a candidate of degree one: $q_c = (y-1)(b_0+b_1x+b_1y+b_1y')$. This simplifies a lot the determination of $\sigma$ because in addition to reducing the number of indeterminate coefficients, it also reduces the number of nonlinear terms. 

In the following section we will treat one case belonging to this one, namely, the case $q=u\,N$ where $u \in \{x,y,y'\}$.

\subsubsection{$q | N$ or $q=u\,N$ where $u \in \{x,y,y'\}$}
\label{qeqn}

In this section we will discuss two very common cases: $q | N$ and $q=y'N$ (or $q=x\,N$ or $q=y\,N$). The main point is that, in addition to appearing frequently (and greatly simplify the calculations) they turn the algorithm $ASymm$ into a full algorithm.

\medskip

\noindent
{\bf The case $q | N$:}

\noindent
To prove the result expressed in the following theorem we only have to consider the case $q = N$ where $p$ and $q$ are not necessarily coprime.

\begin{teor}
\label{fulalg}
Let $\,y''=\phi(x,y,y')=M(x,y,y')/N(x,y,y'),$ where $M$ and $N$ are coprime polynomials, be a rational 2ODE presenting a global analytical first integral $I(x,y,y')$ and a symmetry vector field ${X_e}^{(1)}=\nu\,\partial_y+D_x[\nu]\,\partial_{y'}$ such that the function $\sigma$ defined as $\sigma \equiv -D_x[\nu]/\nu$ is a rational function given by $\sigma={\frac{p}{N}}$, where $p$ is a polynomial. Let $deg_M$, $deg_N$ and $deg_p$ be the degrees of polynomials $M$, $N$ and $p$, respectively. If $deg_M \leq deg_N + 1$ then $deg_p \leq deg_N$. On the contrary, if $deg_M > deg_N + 1$ then $deg_p \leq deg_M - 1.$
\end{teor}

{\bf Proof of the Theorem \ref{fulalg}:} The proof follows directly from the proof of the theorem \ref{relmusigrat}. 

\begin{cor}
\label{fulalgcor}
If the assumptions of the theorem \ref{fulalg} are fulfilled then $ASymm$ becomes a full algorithm.
\end{cor}

{\bf Proof of the Corollary \ref{fulalgcor}:} From the theorem \ref{fulalg} we have that the degree of the polynomial $p$ depends only on the degrees of the polynomials $M$ and $N$. Since they are finite, the maximum degree for the polynomial $p$ is finite and, so, the procedure $ASymm$ will always end.

\medskip

\noindent
{\bf The case $q = y'N$:}

\noindent
This case is very similar to the one previously discussed. To avoid repetition let's just cite the result: 

\medskip

\noindent
{\it Since $q$ is known (i.e., well defined), the maximum degree of $p$ will be limited (as in the previous case) and, as in the previous case, {\em ASymm} becomes a full algorithm.}

\medskip

\noindent
For example, the 2ODE
\begin{equation}
\label{2odeexqequalN}
y'' = -{\frac {{x}^{2}yy'-{x}^{2}{y'}^{2}-x{y}^{3}-x{y}^{2}y'-x{y'}^{2}+{y}^{3}+
{y}^{2}y'+2\,y{y'}^{2}-{y'}^{2}}{y \left( {x}^{2}-y \right) }}
\end{equation}
has a $\sigma$ given by
\begin{equation}
\label{sigexqeqN1}
\sigma = -{\frac {{x}^{2}y'+x{y}^{2}-2\,yy'}{y \left( {x}^{2}-y \right) }},
\end{equation}
where the fact that $q=N$ greatly simplifies the computation of $\sigma$ (The calculation takes 0.5 seconds). With the symmetry 
\begin{equation}
{X_e}^{(1)} = {\rm e}^{-\int_x\left[-{\frac {{x}^{2}y'+x{y}^{2}-2\,yy'}{y \left( {x}^{2}-y \right) }}\right]}\,\left(\partial_y + {\frac {{x}^{2}y'+x{y}^{2}-2\,yy'}{y \left( {x}^{2}-y \right) }}\,\partial_{y'}\right),
\label{xeexperf1}
\end{equation}
we can obtain the first integral
\begin{equation}
{\frac { \left( xy'-{y}^{2} \right) {{\rm e}^{x}}}{xy-y'}}
\label{fiperfex1}
\end{equation}
and the general solution
\begin{equation}
y \left( x \right) =\displaystyle{\frac{{{\rm e}^{\int \!{\frac {C_1\,x}{{{\rm e}^{
x}}x+C_1}}{dx}}}}{\int \!-\displaystyle{\frac{{{\rm e}^{\int \!{\frac {C_1\,x}{{{\rm e}^{x}}x+C_1}}{dx}}}{{\rm e}^{x}}}{{
{\rm e}^{x}}x+C_1}}{dx}+C_2}}
\end{equation}
in less than 0.1 second.

\medskip
\noindent
As another example, the 2ODE
\begin{equation}
\label{2odeqeqzN}
y'' = -{\frac {2\,{x}^{3}{y'}^{5}-4\,{x}^{2}y{y'}^{3}-{x}^{2}{y'}^{4}-{x}^
{2}{y'}^{3}+{x}^{2}{y'}^{2}+2\,x{y}^{2}y'+2\,xy{y'}^{2}-{y}^{2}+yy'-y}{2\,x
 \left( {x}^{3}{y'}^{4}-2\,{x}^{2}y{y'}^{2}-{x}^{2}{y'}^{2}+x{y}^{2}+y
 \right) }}
\end{equation}
has a $\sigma$ given by
\begin{equation}
\label{sigexqeqNz}
\sigma = -{\frac {{x}^{2}{y'}^{4}-{x}^{2}{y'}^{2}-2\,xy{y'}^{2}+{y}^{2}+y}{2xy'
 \left( {x}^{3}{y'}^{4}-2\,{x}^{2}y{y'}^{2}-{x}^{2}{y'}^{2}+x{y}^{2}+y
 \right) }},
\end{equation}
that can be calculated in 4.3 seconds. From the knowledge of $\sigma$ we can calculate the first integral
$$
{\frac {\ln  \left( {x}^{2}{y'}^{2}-y \right) x{y'}^{2}-y\ln  \left( {x}
^{2}{y'}^{2}-y \right) +1}{x{y'}^{2}-y}}
$$
in about one second.

\medskip

\subsection{$M$ and $p$ present common factors or common monomials}
\label{meqp}

Let's return to the equation (\ref{eqpq4}): the term 
$$
\frac{p^2N + p\,q\,M_{y'} - q^2M_{y} - q\,p_xN - y'q\,p_yN + p\,q_xN + y'p\,q_yN}{M}
$$
is a polynomial. We can write it as
\begin{equation}
\label{ovm2}
p\,\frac{p\,N + q\,M_{y'}+ q_xN + y'q_yN}{M} - q\,\frac{q\,M_{y} + p_xN + y'p_yN}{M},
\end{equation}

\noindent
and, in this form we can see that in many cases some factors of $M$ and $p$ may be common. Besides, since $\phi$ and $\sigma$ are given, respectively, by 
\begin{equation}
\phi = - \frac{\partial_x[I]+z\,\partial_y[I]}{\partial_{y'}[I]} \,\,\, {\rm and} \,\,\, \sigma = \frac{\partial_y[I]}{\partial_{y'}[I]} \,\, \Rightarrow \,\,
\phi = - \frac{\partial_x[I]}{\partial_{y'}[I]} + y'\frac{p}{q}
\end{equation}
it is very likely (in fact, this is very often the case) that the polynomials $M$ and $p$ have many common monomials (case $N \neq q$) or that the polynomials $M_z$ and $p$ have many common monomials (case $N = q$). 

Each of these situations greatly simplifies the problem. Let's take, for example, the nonlinear 2ODEs 94 and 183 (see section \ref{perfitK}) of Kamke's book:

\bigskip

\noindent
2ODE 94:

$$
y''={\frac {2\,{y}^{3}{x}^{3}-2\,yy'{x}^{3}-3\,{x}^{2}{y}^{2}-axy-9\,{x}^{2}y'-b}{2{x}^{3}}},
$$

$$
\sigma={\frac {-2\,{y}^{3}{x}^{3}+2\,yy'{x}^{3}+3\,{x}^{2}{y}^{2}+axy+5\,{x}^{2}y'+b}{2{x}^{2} \left( xy'+y \right) }}.
$$

\bigskip

\noindent
2ODE 183:

$$
y''={\frac {{x}^{2}{y'}^{2}+{x}^{2}-{y}^{2}}{2{x}^{2}y}},\,\,\,\,\,\,\sigma={\frac {{x}^{2}{y'}^{2}+{x}^{2}-{y}^{2}}{-2xy \left( xy'-y \right) }}.
$$

\bigskip

\noindent
We can see that in the 2ODE 94, the monomials in the polynomias $M$ and $p$ are the same, while in 2ODE 183 we have that $M=p$. Also in the two 2ODEs $N$ and $q$ have common factors ($2x^2$ in 2ODE 94 and $2xy$ in 2ODE 183). Therefore, the time taken to calculate $\sigma$ is (in both 2ODEs) tenths of a second. 

\bigskip

As another example, for the 2ODE (\ref{2odeexqequalN}) ($q=N$) 
$$
y'' = -{\frac {{x}^{2}yy'-{x}^{2}{y'}^{2}-x{y}^{3}-x{y}^{2}y'-x{y'}^{2}+{y}^{3}+
{y}^{2}y'+2\,y{y'}^{2}-{y'}^{2}}{y \left( {x}^{2}-y \right) }},
$$
using the monomials present in $M_{y'}=-{x}^{2}y+2\,{x}^{2}y'+x{y}^{2}+2\,xy'-{y}^{2}-4\,y'y+2\,y'$ we can calculate $p=- {x}^{2}y'-x{y}^{2}+2\,yy'$ (and so, $\sigma$) almost instantly.

\newpage

\section{Performance}
\label{perfor}

In this section we present the performance of our procedure in two types of arena. The first one consists of a well-known set of 2ODEs, namely the nonlinear 2ODEs of E. Kamke's well-known Handbook of exact solutions for ODEs \cite{Kam}. The second test arena consists of a set of 2ODEs that mostly exhibit nonlocal symmetries.

\medskip

\subsection{Nonlinear rational 2ODEs of Kamke's Handbook}
\label{perfitK}

In this section we show our algorithm dealing with the differential equations of Kamke's handbook: by this we mean the nonlinear rational 2ODEs that present (at least) a Liouvillian first integral. Also, let's just consider the subset of those that don't have trivial symmetries like $[\xi=1, \eta=0]$ or $[\xi=0, \eta=1]$, i.e. let's leave out 2ODEs so that $\phi$ doesn't explicitly present one of the variables $(x,y)$. Let's refer to them by the number given in Kamke's book. They are the 2ODEs referenced by the numbers:
$\!\![78, 79, 80, 87, 90, 92, 93, 94, 97, 98, 99,\! 108,\! 133,\! 156,\! 169,\! 172,\! 173,\! 174, 175, 176,$ $178,179, 180, 181, \,182, \,183, \,184, \,185,$ 
$189, 190, 193, 206, 226, 227, 228, 229, 231]$.
\begin{table}[h]
{\begin{center} {\footnotesize
\begin{tabular}
{|c|c|c|c|}
\hline
ODE & $\phi$ & $\sigma$ & $\xi,\,\eta,\,\nu$ \\	
\hline
78 & $-{\frac { \left( y-1 \right) y'}{x}}$ & ${\frac {y-2}{x}}$ & $x,\,\,0,\,\,-xy'$  \\
\hline
79 & $-{\frac {-{x}^{2}{y'}^{2}+{y}^{2}+2\,y'}{x}}$ & $-xy'+y$ & $x,\,\,-y,\,\,-xy'-y$  \\
\hline
80 & $-{\frac {a{x}^{2}{y'}^{2}-2\,yy'ax+a{y}^{2}-b}{x}}$ & $-{\frac {1}{x}}$ & $0,\,\,x,\,\,x$  \\
\hline
87 & $-{\frac {ay{y'}^{2}+bx}{{x}^{2}}}$ & ${\frac {ay{y'}^{2}+bx}{x \left( xy'-y \right) }}$ & $x,\,\,y,\,\,-xy'+y$  \\
\hline
90 & ${\frac {{x}^{4}{y'}^{2}-4\,y}{4{x}^{2}}}$ & ${\frac {{x}^{4}{y'}^{2}+12\,xy'-4\,y}{-4x \left( xy'+2\,y \right) }}$ & $x,\,\,-2y,\,\,-xy'-2y$  \\
\hline
92 & ${\frac {{y}^{3}{x}^{3}-yy'{x}^{3}-12\,yx-24}{{x}^{3}}}$ & $-{\frac {{y}^{3}{x}^{3}-yy'{x}^{3}+2\,{x}^{2}y'-12\,yx-24}{{x}^{2}
 \left( xy'+y \right) }}
$ & $-x,\,\,y,\,\,xy'+y$  \\
\hline
93 & ${\frac {a \left( -xy'+y \right) ^{2}}{{x}^{3}}}$ & $-{\frac {1}{x}}$ & $0,\,\,x,\,\,x$  \\
\hline
94 & ${\frac {2\,{y}^{3}{x}^{3}-2\,yy'{x}^{3}-3\,{x}^{2}{y}^{2}-axy-9\,{x}^{2}y'-b}{2{x}^{3}}}$ & ${\frac {-2\,{y}^{3}{x}^{3}+2\,yy'{x}^{3}+3\,{x}^{2}{y}^{2}+axy+5\,{x}^{2}y'+b}{2{x}^{2} \left( xy'+y \right) }}$ & $-x,\,\,y,\,\,xy'+y$  \\
\hline
97 & $-{\frac {-y'{x}^{3}-2\,xyy'+4\,{y}^{2}}{{x}^{4}}}$ & $-{\frac {2y}{x^3}}$ & $x,\,\,2y,\,\,-xy'+2y$  \\
\hline
98 & $-{\frac {-y'{x}^{3}-{x}^{2}{y'}^{2}+4\,{y}^{2}}{{x}^{4}}}$ & $-{\frac {xy'+2\,y}{{x}^{3}}}$ & $x,\,\,y,\,\,y-xy'/2$  \\
\hline
99 & ${\frac { \left( -xy'+y \right) ^{3}}{{x}^{4}}}$ & $-{\frac {1}{x}}$ & $0,\,\,x,\,\,x$  \\
\hline
108 & $-{\frac {-ax+{y'}^{2}-b}{y}}$ & ${\frac {2\,{a}^{2}{x}^{2}-2\,ax{y'}^{2}+4\,bxa-ayy'-2\,b{y'}^{2}+2\,{b}^
{2}}{-y \left( 2\,xy'a-3\,ay+2\,by' \right) }}$ & $\frac{2(ax+b)}{3a},\,y,\,{\frac {2\,xy'a-3\,ay+2\,by'}{3a}}$  \\
\hline
133 & $-{\frac {y' \left( y'-1 \right) }{x+y}}$ & ${\frac {y' \left( y'-1 \right) }{ \left( x+y \right)  \left( 1+y'
 \right) }}$ & $1,\,\,-1,\,\,-1-y'$  \\
\hline
156 & ${\frac {a{x}^{2}+bx+2\,{y'}^{2}+c}{3y}}$ & see below & see below  \\
\hline
169 & $-{\frac {y' \left( xy'-y \right) }{yx}}$ & $-{\frac {y'}{y}}$ & $0,\,\,y,\,\,y$  \\
\hline
172 & ${\frac {-bx{y}^{3}-ayy'+x{y'}^{2}}{yx}}$ & ${\frac {bx{y}^{3}+ayy'-x{y'}^{2}-3\,yy'}{y \left( xy'+2\,y \right) }}$ & $x,\,\,-2y,\,\,-xy'-2y$  \\
\hline
173 & $-{\frac {y' \left( ay+2\,xy' \right) }{yx}}$ & $-\frac {y'}{y}$ & $0,\,\,y,\,\,y$  \\
\hline
174 & ${\frac {y' \left( 2\,xy'-y-1 \right) }{yx}}$ & $-{\frac {2\,xy'-1}{yx}}$ & $x,\,\,0,\,\,-xy'$  \\
\hline
175 & ${\frac {y' \left( -ay+2\,xy' \right) }{yx}}$ & $-{\frac {y'}{y}}$ & $0,\,\,y,\,\,y$  \\
\hline
176 & ${\frac {4y' \left( xy'-y \right) }{yx}}$ & $-{\frac {y'}{y}}$ & $0,\,\,y,\,\,y$  \\
\hline
178 & $-{\frac {x{y'}^{2}+xy'-yy'-y}{x \left( x+y \right) }}$ & ${\frac {xy'-x-2\,y}{x \left( x+y \right) }}$ & $-x,\,\,x,\,\,xy'+x$  \\
\hline
\end{tabular} }
\end{center}}
\caption{2ODEs, functions $\sigma$ and symmetries}
\label{pss}
\end{table}

\begin{table}[h]
{\begin{center} {\footnotesize
\begin{tabular}
{|c|c|c|c|}
\hline
ODE & $\phi$ & $\sigma$ & $\xi,\,\eta,\,\nu$ \\	
\hline
179 & ${\frac {y' \left( xy'-y \right) }{2yx}}$ & $-{\frac {y'}{y}}$ & $0,\,\,y,\,\,y$  \\
\hline
180 & ${\frac {2\,({x}^{2}{y'}^{2}+yy'x-{y}^{3}-xy'+2\,{y}^{2}-y)}{{x}^{2} \left( y-1 \right) }}$ & $-{\frac {y' \left( 2\,y-1 \right) }{y \left( y-1 \right) }}$ & $0,\,y^2-y,\,\,y^2-y$  \\
\hline
181 & ${\frac {{x}^{2}{y'}^{2}-2\,yy'x+{y}^{2}}{{x}^{2} \left( x+y \right) }}$ & $-{\frac {xy'-y}{x \left( x+y \right) }}$ & $x,\,\,y,\,\,y-xy'$  \\
\hline
182 & ${\frac {a \left( {x}^{2}{y'}^{2}-2\,yy'x+{y}^{2} \right) }{{x}^{2} \left( y-x \right) }}$ & ${\frac { \left( xy'-y \right) a}{x \left( x-y \right) }}$ & $x,\,\,y,\,\,y-xy'$  \\
\hline
183 & ${\frac {{x}^{2}{y'}^{2}+{x}^{2}-{y}^{2}}{2{x}^{2}y}}$ & ${\frac {{x}^{2}{y'}^{2}+{x}^{2}-{y}^{2}}{-2xy \left( xy'-y \right) }}$ & $x,\,\,y,\,\,y-xy'$  \\
\hline
184 & $-{\frac {b{x}^{2}{y'}^{2}+cxyy'+d{y}^{2}}{a{x}^{2}y}}$ & $-{\frac {y'}{y}}$ 	& $0,\,\,y,\,\,y$  \\
\hline
185 & ${\frac {{y'}^{2}{x}^{3}+{y}^{2}ax-2\,{x}^{2}yy'+2\,{x}^{2}{y'}^{2}+2\,a{y}^{2}-4\,yy'x+x{y'}^{2}-2\,yy'}{yx \left( {x}^{2}+2\,x+1 \right) }}$ & $-{\frac {y'}{y}}$ & $0,\,\,y,\,\,y$  \\
\hline
189 & $-{\frac {y{y'}^{2}+ax}{{y}^{2}}}$ & ${\frac { \left( y{y'}^{2}+ax \right) x}{{y}^{2} \left( xy'-y \right) }}$ & $x,\,\,y,\,\,y-xy'$  \\
\hline
190 & $-{\frac {y{y'}^{2}-ax-b}{{y}^{2}}}$ & $-{\frac { \left( -y{y'}^{2}+ax+b \right)  \left( ax+b \right) }{{y}^{2} \left( xy'a-ay+by' \right) }}$ & $x\!+\!\frac{a}{b},\,y,\,\!{\frac {ay-xy'a-by'}{a}}$  \\
\hline
193 & ${\frac {y' \left( 2\,{y}^{2}{y'}^{2}-2\,x{y'}^{2}+4\,yy'+1 \right) }{-({y}^{2}+x)}}$ & ${\frac { 2\left( {y}^{2}y'-xy'+y \right) y'}{{y}^{2}+x}}$ & $-2y,\,1,\,2yy'+1$  \\
\hline
206 & ${\frac {y' \left( y'y{a}^{2}-{x}^{2}yy'-{a}^{2}x+x{y}^{2} \right) }{-{a}^{4}+{a}^{2}{x}^{2}+{a}^{2}{y}^{2}-{x}^{2}{y}^{2}}}$ & ${\frac {yy'}{{a}^{2}-{y}^{2}}}$ & \scriptsize{$0$},\scriptsize{$\sqrt{a^2-y^2}$},\scriptsize{$\sqrt{a^2-y^2}$}  \\
\hline
226 & ${\frac {yx \left( xy'+y \right) }{y'}}$ & $-{\frac {y'}{y}}$ & $0,\,\,y,\,\,y$  \\
\hline
227 & $-{\frac {4\,{y'}^{2}}{xy'-y}}$ & $-{\frac {y'}{y}}$ & $0,\,\,y,\,\,y$  \\
\hline
228 & ${\frac {{y'}^{4}+2\,{y'}^{2}+1}{xy'-y}}$ & $-{\frac { \left( {y'}^{2}+1 \right) y'}{xy'-y}}$ & $-y,\,\,x,\,\,yy'+x$  \\
\hline
229 & $-{\frac {b{y}^{2}}{ay'{x}^{3}}}$ & $-{\frac {y'}{y}}$ & $0,\,\,y,\,\,y$  \\
\hline
231 & $-{\frac {2\,y{y'}^{3}+3\,xy'+y}{2\,{y}^{2}y'+{x}^{2}}}$ & ${\frac { \left( 2\,y{y'}^{3}+3\,xy'+y \right) x}{ \left( 2\,{y}^{2}y'+{x}
^{2} \right)  \left( xy'-y \right) }}$ & $x,\,y,\,y-xy'$  \\
\hline
\end{tabular} }
\end{center}}
\caption{2ODEs, functions $\sigma$ and symmetries  (Table 1: cont.)}
\label{pss2}
\end{table}

\noindent
$\sigma_{156}={\frac {2\,{a}^{2}{x}^{4}+4\,ab{x}^{3}+4\,a{x}^{2}{y'}^{2}+4\,ac{x}^{2}-6\,axyy'+2\,{b}^{2}{x}^{2}+4\,bx{y'}^{2}-18\,a{y}^{2}+4\,bcx-3\,byy'+4\,c{y'}^{2}+2\,{c}^{2}}{-3y \left( 2\,a{x}^{2}y'-6\,axy+2\,xy'b-3\,by+2\,cy' \right) }}$

\medskip

\noindent
$\xi_{156}={\frac {a{x}^{2}+bx+c}{b}},\,\,\,\eta_{156}=3y\,\left({\frac {ax}{b}}+\frac{1}{2}\right),\,\,\,\nu_{156}={\frac {2\,a{x}^{2}y'-6\,axy+2\,xy'b-3\,by+2\,cy'}{-2b}}$

\bigskip

\noindent
Some comments:

\begin{itemize}
\item Our method was able to find the symmetries for all nonlinear rational 2ODEs in Kamke's book that presented (at least) a Liouvillian first integral. 

\item There are some nonlinear 2ODEs of Kamke's book that present generic functions of $x$, $y$ or $y'$. In order to test our method we have made several substitutions of these generic functions for rational functions. We realize that the algorithm is able to handle well all 2ODEs that had a Liouvillian first integral. 

\item The 37 2ODEs (that fit the conditions for which the algorithm was built) have very simple symmetries. The same result was found for the 2ODEs presenting arbitrary functions (when replaced by rational functions).

\item The time spent calculating the $\sigma$ function was not put in the table because it was too short (tenths of a second in most cases). Only at 2ODEs 108 and 156 did the algorithm spend a few seconds.

\item The time spent calculating the symmetries by the Maple Computer Algebra System (CAS) {\tt symgen} command (see \cite{Noscpc1998}) was even shorter.
\end{itemize}			

Although satisfied with the result presented by the algorithm in this test arena, we conclude that it is not suitable for testing the method performance in dealing with 2ODEs that have nonlocal symmetries.

\medskip

\subsection{Rational 2ODEs presenting nonlocal symmetries}
\label{perfitO}

In this section we will look at the performance of our algorithm in a test arena made up of a set of rational 2ODEs that do not have local symmetries. For the first two 2ODEs we will present the steps of the algorithm in more detail. 

\begin{itemize}
\item {\bf 2ODE$_1$:} 
\begin{equation}
\label{2odesnl1}
y''=-{\frac {xyy'-2\,x{y'}^{2}+yy'-{y'}^{2}-y+2\,y'}{xy-1}}
\end{equation}

\noindent
so, the steps of the algorithm are:

\begin{enumerate}
\item Starting with the hypothesis that $q$ divides $N$ (which turns $ASymm$ into a full algorithm and greatly simplifies the calculations) we have as a consequence of the theorem \ref{fulalg} that the maximum degree of $p$ is 2 and therefore we will have as candidates:

$p_c = a_4x^2+a_5y^2+a_6y'^2+a_7xy+a_8yy'+a_9xy'+a_1x+a_2y+a_3y'+a_0$,

$q_c = N = xy-1$.

\item Substituting $\sigma$ for $p_c/q_c$ in equation (\ref{eqproofpq2}), collecting the resulting polynomial equation in $(x,y,y')$ and equalizing the coefficients of each monomial to zero we get the following system $AE$ of algebraic equations for $a_i$:
$$
AE = \{a_4^2, a_5^2, a_6^2, 2 a_1 a_4, 2 a_2 a_5, 2 a_3 a_6, 2 a_5 a_9, 2 a_6 a_9, 2 a_4 a_7+5 a_4, 2 a_4 a_8-a_4,
$$
$$
2 a_5 a_8-a_5,\! 2 a_6 a_7+a_6, 2 a_5 a_6+a_9^2,\! 2 a_0 a_4+a_1^2-2 a_4,\! 2 a_0 a_1-2 a_1+2 a_4, 2 a_4 a_5+a_8^2 
$$
$$
-a_8,\!2 a_4 a_9+2 a_7 a_8+4 a_8, 2 a_5 a_7+2 a_8 a_9+3 a_5, a_0^2-2 a_0+a_1+1, 2 a_4 a_6+a_7^2+3 a_7
$$
$$
+2,2 a_0 a_2-2 a_2-a_3+a_8, 2 a_0 a_5+a_2^2-2 a_5-a_9,\! 2 a_1 a_8+2 a_2 a_4-a_1-2 a_4, \!2 a_1 a_5+
$$
$$
2 a_2 a_8-a_2-a_8,\, 2 a_1 a_7+\,2 a_3 a_4+5 a_1+2 a_4,\, 2 a_2 a_9+2 a_3 a_5+2 a_5+a_9,\, 2 a_2 a_6\,+
$$
$$
2 a_3 a_9+2 a_6+a_9, 2 a_6 a_8+2 a_7 a_9+a_6+2 a_9, 2 a_0 a_3+2 a_0+a_2+a_7-1, 2 a_0 a_6+a_3^2+
$$
$$
\!a_3+2 a_6+\!a_9, \!2 a_1 a_6+2 a_3 a_7+3 a_3+a_7+\!1\!,\! 2 a_1 a_9+2 a_2 a_7+2 a_3 a_8+4 a_2+2 a_8,\!2 a_0 a_8
$$
$$
+2 a_1 a_2-a_0-a_1\!-a_7-2 a_8,\! 2 a_0 a_7+2 a_1 a_3+5 a_0+2 a_1+a_8-3,\! 2 a_0 a_9+2 a_2 a_3+2 a_2
$$
$$
+a_3+2 a_5-2 a_6\}.
$$

\item The $AE$ system has as a solution
$$
\{a_0 \!=\! 1, a_1 \!=\! 0, a_2 \!=\! 0, a_3 \!=\! 0, a_4 \!=\! 0, a_5 \!=\! 0, a_6 \!=\! 0, a_7 \!=\! -1, a_8 \!=\! 0, a_9 \!=\! 0\}.
$$

\item Substituting in $p_c/q_c$ we have that $\sigma=\displaystyle{-{\frac {xy'-1}{xy-1}}}$.

\end{enumerate}

Therefore, the 2ODE (\ref{2odesnl1}) admits the symmetry 

\begin{equation}
{X_e}^{(1)} = {\rm e}^{-\int_x\left[-{\frac {xy'-1}{xy-1}}\right]}\,\left(\partial_y + {\frac {xy'-1}{xy-1}}\,\partial_{y'}\right).
\label{xeexnl1}
\end{equation}
Using the Lie method we can determine the first integral 
\begin{equation}
{\frac { \left( y-y' \right) {{\rm e}^{-x}}}{xy'-1}}
\label{fiexnl1}
\end{equation}
and the solution
\begin{equation}
y \left( x \right) = \left( \int \!\frac{{C_1}\,{{\rm e}^{-\int \!{
\frac {{{\rm e}^{-x}}}{{C_1}\,x+{{\rm e}^{-x}}}}{dx}}}}{ \left( {C_1}\,x+{{\rm e}^{-x}} \right)}{dx}+{C_2} \right) {
{\rm e}^{\int \!{\frac {{{\rm e}^{-x}}}{{C_1}\,x+{{\rm e}^{-x}}}}
{dx}}}
\label{solexnl1}
\end{equation}

\item{\bf 2ODE$_2$:} 
\begin{equation}
\label{2odesnl2}
y''={\frac {{x}^{2}{y}^{2}+{x}^{2}yy'-2\,{y'}^{2}xy-x{y'}^{3}+{y'}^{4}-{x
}^{2}y'+x{y}^{2}-y{y'}^{2}-yx}{2y' \left( yx-{y'}^{2}-x \right) }}
\end{equation}

\noindent
so, the steps of the algorithm are:

\begin{enumerate}
\item Starting with the hypothesis that $q$ divides $N$ (and that it is a factor of $N$) we will get a positive solution with the candidates:

$p_c = a_1x+a_2y+a_3y'+a_0,$

$q_c = y'$.

\item Substituting $\sigma$ for $p_c/q_c$ in equation (\ref{eqproofpq2}), collecting the resulting polynomial equation in $(x,y,y')$ and equalizing the coefficients of each monomial to zero we get the following system $AE$ of algebraic equations for $a_i$:
$$
AE = \{-4 a_0^2, -2 a_0^2, -2 a_2^2, 4 a_2^2, -4 a_1 a_3, -4 a_1^2-2 a_1, -2 a_1^2-a_1, 4 a_1^2+2 a_1, 
$$
$$
-4 a_0 a_1-a_0, -4 a_2 a_3+2 a_2, 8 a_2 a_3-4 a_2, -4 a_0 a_2+a_3, -8 a_1 a_3+4 a_1+2, -4 a_1 a_3
$$
$$
+2 a_1+1, 8 a_1 a_3-4 a_1-2, -2 a_3^2+2 a_2+a_3, -4 a_0 a_3+2 a_1+1, -4 a_0 a_2+4 a_2^2+a_3, 
$$
$$
8 a_0 a_2-4 a_2^2-2 a_3, 4 a_0^2-8 a_0 a_2+2 a_3, 4 a_0^2-4 a_0 a_2+a_3, -4 a_1 a_2-a_2+a_3, 8 a_0 a_1
$$
$$
-4 a_1 a_2+2 a_0-a_2, -8 a_0 a_3+4 a_0+4 a_1+2, -4 a_0 a_3+2 a_0+2 a_1+1, -2 a_0^2+8 a_0 a_2
$$
$$
-2 a_2^2-2 a_3, -4 a_1 a_2+4 a_3^2-5 a_2-a_3, 8 a_1 a_2-2 a_3^2+4 a_2-a_3, -8 a_0 a_1-2 a_3^2-2 a_0
$$
$$
+2 a_2,\! 8 a_0 a_3-4 a_2 a_3-4 a_0-4 a_1-2,\!\! -4 a_0 a_1-4 a_3^2-a_0+4 a_2+3 a_3,\!\! -4 a_0 a_1+8 a_1 a_2
$$
$$
-a_0+2 a_2-a_3, -4 a_0 a_3+8 a_2 a_3+2 a_0+2 a_1-4 a_2+1, 8 a_0 a_3-8 a_2 a_3-4 a_0-4 a_1+
$$
$$
4 a_2-2, 8 a_0 a_1-8 a_1 a_2+4 a_3^2+2 a_0-6 a_2-2 a_3\}.
$$

\item The solution of the $AE$ system is
$$
\{a_0 = 0, a_1 = -1/2, a_2 = 0, a_3 = 0\}.
$$

\item Substituting in $p_c/q_c$ we have that $\sigma=\displaystyle{-{\frac {x}{2y'}}}$.

\end{enumerate}

Therefore, the 2ODE (\ref{2odesnl2}) admits the symmetry 

\begin{equation}
{X_e}^{(1)} = {\rm e}^{-\int_x\left[-{\frac {x}{2y'}}\right]}\,\left(\partial_y + {\frac {x}{2y'}}\,\partial_{y'}\right),
\label{xeexnl2}
\end{equation}
leading to the the first integral 
\begin{equation}
{\it Ei} \left( 1,- \frac{1}{ yx-{y'}^{2}} \right) +x\,{{\rm e}
^{ \left({ \frac{1}{ yx-{y'}^{2}}} \right) }}.
\label{fiexnl2}
\end{equation}

\end{itemize}

In what follows we will present a table with the 2ODE ($\phi$), the function $\sigma$ together with the time spent by $ASymm$ to determine it.

\begin{table}[h]
{\begin{center} {\footnotesize
\begin{tabular}
{|c|c|c|c|}
\hline
2ODE & $\phi$ & $\sigma$ & Time \\	
\hline
1 & $-{\frac {{x}^{3}{y'}^{3}-{x}^{2}y{y'}^{2}-x{y}^{2}y'+{y}^{3}}{-{x}^{2}{y'}^{2}+{x}^{2}y+2\,yy'x-xy'-{y}^{2}}}$ & 
${\frac {{x}^{3}{y'}^{2}-2\,{x}^{2}yy'+x{y}^{2}-yx+y'}{-{x}^{2}{y'}^{2}+{x}^{2}y+2\,yy'x-xy'-{y}^{2}}}$ & 0.95  \\
\hline
2 & ${\frac {{x}^{2}{y'}^{2}-2\,x{y'}^{3}+{y'}^{4}-{y}^{2}+y'y}{y \left( {x}^{2}-2\,y'x+{y'}^{2}-y+y' \right) }}$ & 
$-{\frac {y' \left( {x}^{2}-2\,y'x+{y'}^{2} \right) }{y \left( {x}^{2}-2\,y'x+{y'}^{2}-y+y' \right) }}$ & 0.2  \\
\hline
3 & ${\frac {{x}^{2}{y'}^{2}-2\,x{y}^{2}y'-2\,xy{y'}^{2}+{y}^{4}+2\,{y}^{3}y'+2\,xyy'+x{y'}^{2}-{y}^{2}y'-y'x}{-x \left( y'x-{y}^{2}-x \right) }}$ & 
$-{\frac {2y}{x}}$ & 0.17  \\
\hline
4 & ${\frac { \left( -y'y+x+y' \right)  \left( {y'}^{2}-1 \right) }{{y}^{2}{y'}^{2}-2\,xyy'+{y}^{2}y'+{x}^{2}-xy-y'y}}$ & 
${\frac {-y' \left( -yy'+x+y' \right) }{{y}^{2}{y'}^{2}-2\,xyy'+{y}^{2}y'+{x}^{2}-xy-yy'}}$ & 0.07  \\
\hline
5 & ${\frac { \left( xy-{y'}^{2}+y' \right)  \left( xy'+y \right) }{2\,xyy'-2\,{y'}^{3}+xy+{y'}^{2}}}$ & 
$-{\frac {x \left( xy-{y'}^{2}+y' \right) }{2\,xyy'-2\,{y'}^{3}+xy+{y'}^{2}}}$ & 0.09  \\
\hline
\end{tabular} }
\end{center}}
\caption{2ODEs with non local symmetries and computation of the $\sigma$'s}
\label{nostab}
\end{table}

\bigskip

\noindent
A few more comments:

\begin{itemize}
\item In this arena, our method was able to find the symmetries for rational 2ODEs presenting only non local symmetries in a very short time. 

\item Neither of these 2ODEs can be solved/reduced by the dsolve command (the Maple ODE solver).

\item Our approach works `better' (i.e., it comparatively spents shorter times ) when the symmetries are non local.
\end{itemize}

\medskip

\section{Applications}
\label{applic}

Nonlinear 2ODEs model a vast amount of physical, chemical and biological phenomena, directly or as part of the final process of solving systems of partial differential equations (PDE systems) that mathematically model more complex phenomena. Within the set of nonlinear 2ODEs, a category that occupies a great highlight consists of nonlinear oscillators, usually connected to a source that provides the system with periodic excitation. These nonlinear oscillators generally present several parameters that represent magnitudes involved in the phenomenon described and, in general, present chaotic behavior for a certain region of the space of the parameters. Therefore, it is very important to determine the values of the parameters for which the system is integrable (values around which it is possible to make a progressive variation of the parameters to study the bifurcation of the trajectories and the transition from the regular to the chaotic regime). One of the great advantages of our method is that it can perform an integrability analysis (when the 2ODE depends on parameters) if we consider that some algebraic combination of them can lead to new solutions. In the next subsections we will present some applications of our method to some 2ODEs representing nonlinear oscillators.

\subsection{Helmholtz oscillator with friction}
\label{helmowf}

The 2ODE that represents the Helmholtz oscillator with friction \cite{AlSa} is 
\begin{equation}
y'' = a\,y'+b\,y-c\,y^2,
\label{2odehowf}
\end{equation}
where $a$, $b$ and $c$ are parameters. Since the 2ODE (\ref{2odehowf}) does not depend explicity on $x$ we have $\partial_x$ as a symmetry. This permits a order reduction as shown in section \ref{ops}. In the analysis that follows, we will skip this more trivial case and focus our attention on possible relationships between the parameters of the 2ODE that allow us to obtain Liouvillian first integrals.

\medskip

For the choice 
$$
p_c = a_4x^2+a_5y^2+a_6y'^2+a_7xy'+a_8xy+a_9yy'+a_1x+a_2y+a_3y'+a_0,
$$
$$
q_c = b_1x+b_2y+b_3y'+b_0,
$$
solving the $AE$ system for the coefficients of $p_c$ and $q_c$ and for the parameters $a$, $b$ and $c$ we obtain three solutions (we're just considering the solutions where $c$ is not null, i.e., where the 2ODE (\ref{2odehowf}) is not linear):

\medskip

\noindent
{\bf Case 1:} $b={\frac{6}{25}}\,a^2$:

\medskip

\noindent
For this solution, the 2ODE becomes
\begin{equation}
\label{ode2cas2}
y''=-c{y}^{2}+ay'+{\frac {6}{25}}\,{a}^{2}y
\end{equation}
and
\begin{equation}
\label{scas2}
\sigma={\frac {12\,{a}^{4}-200\,{a}^{2}c\,y+625\,{c}^{2}{y}^{2}-250\,ac\,y'}{5\,(12\,{a}^{3}-50\,ac\,y+125\,c\,y')}}.
\end{equation}
Using the symmetry 
\begin{equation}
{X_e}^{\!(1)}\! =\! {\rm e}^{\!\int_x\!\left[\frac {\!12\,{a}^{4}\!-\!200\,{a}^{2}cy\!+\!625\,{c}^{2}{y}^{2}\!-\!250\,acy'}{-5\,(12\,{a}^{3}-50\,acy+125\,cy')}\!\right]}\,\!\!\left(\!\partial_y\! -\! \frac {12\,{a}^{4}-200\,{a}^{2}cy+625\,{c}^{2}{y}^{2}-250\,acy'}{5\,(12\,{a}^{3}-50\,acy+125\,cy')}\,\partial_{y'}\!\right)\!,
\label{xcas2}
\end{equation}
we can determine the first integral
\begin{equation}
\label{ficas2}
\left( 72\,{a}^{4}y-600\,{a}^{2}c{y}^{2}+1250\,{c}^{2}{y}^{3}+360\,{a}^{3}y'-1500\,acyy'+1875\,c{y'}^{2} \right)\,{{\rm e}^{-\frac{6a}{5}x}}
\end{equation}
leading to the reduced 1ODE
\begin{equation}
\label{oder1cas2}
y'=\frac{2a}{5}y-{\frac {12a^3}{125c}}+{\frac {\sqrt {432\,{a}^{6}-5400\,{a}^{4}cy+22500\,{a}^{2}{c}^{2}{y}^{2}-31250\,{c}^{3}{y}^{3} +75C_1c\,{{\rm e}^{\frac{6a\,x}{5}}}}}{125\sqrt {3}\,c}},
\end{equation}
which admits the Lie symmetry
\begin{equation}
\label{oder1cas2sympoint}
{{\rm e}^{-\frac{a}{5}x}}\,\partial_x-{\frac {2 \left( 6\,{a}^{2}-25\,cy \right) a\,{{\rm e}^{-\frac{a}{5}x}}}{125\,c}}\,\partial_y.
\end{equation}
Using the symmetry (\ref{oder1cas2sympoint}) above we can integrate the reduced 1ODE (\ref{oder1cas2}) in order to obtain the general solution of the 2ODE (\ref{ode2cas2}):
\begin{equation}
\label{gsoder2cas2}
y ={\frac {6\,a^2}{25\,c}}+{\it RootOf}
 \left(\int ^{{\it \_Z}}\!{\frac {\sqrt {-3750\,{c}^{3}{{\it \_a}}^{3}+3\,c\,C_1}}{-1250\,{{\it \_a}}^{3}{c}^{2}+C_1}}{d{
\it \_a}}+C_2-\frac{{\rm e}^{\frac{a}{5}x}}{5\,a} \right) {{\rm e}^{\frac{2a}{5}x}}.
\end{equation}

\bigskip

\noindent
{\bf Case 2:} $b=-{\frac{6}{25}}\,a^2$:

\medskip

\noindent
For this solution, the 2ODE becomes
\begin{equation}
\label{ode2cas3}
y''=-c{y}^{2}+ay'-{\frac {6}{25}}\,{a}^{2}y
\end{equation}
and
\begin{equation}
\label{scas3}
\sigma={\frac {4\,{a}^{2}y+25\,c{y}^{2}-10\,ay'}{-5\,(2\,ay-5\,y')}}.
\end{equation}
Using the symmetry 
\begin{equation}
{X_e}^{\!(1)}\! =\! {\rm e}^{\!\int_x\!\left[{\frac {4\,{a}^{2}y+25\,c{y}^{2}-10\,ay'}{5\,(2\,ay-5\,y')}}\!\right]}\,\!\!\left(\!\partial_y\! +\! {\frac {4\,{a}^{2}y+25\,c{y}^{2}-10\,ay'}{5\,(2\,ay-5\,y')}}\,\partial_{y'}\!\right)\!,
\label{xcas3}
\end{equation}
we can determine the first integral
\begin{equation}
\label{ficas3}
\left( 12\,{y}^{2}{a}^{2}+50\,c{y}^{3}-60\,ayy'+75\,{y'}^{2} \right)\,{{\rm e}^{-\frac{6a}{5}x}}
\end{equation}
and therefore the reduced 1ODE
\begin{equation}
\label{oder1cas3}
y'=\frac{2\,a}{5}y+\sqrt {-\frac{2\,c}{3}{y}^{3}+{\frac {C_1}{75}}\,{{\rm e}^{\frac{6\,a}{5}x}}},
\end{equation}
which admits the Lie symmetry
\begin{equation}
\label{oder1cas3sympoint}
{{\rm e}^{-\frac{a}{5}x}}\,\partial_x+\frac{2\,a}{5}{{\rm e}^{-\frac{a}{5}x}}y\,\partial_y.
\end{equation}
Using the symmetry (\ref{oder1cas3sympoint}) above we can integrate the reduced 1ODE (\ref{oder1cas3}) in order to obtain the general solution of the 2ODE (\ref{ode2cas3}):
\begin{equation}
\label{gsoder2cas3}
y = {{\it RootOf} \left(\int ^{{\it \_Z}}\!{\frac {1}{\sqrt {-150\,c{{\it \_a}}^{3}+3\,C_1}}}{d{\it \_a}}+C_2-\frac{{\rm e}^{\frac{a}{5}x}}{3\,a} \right) }{{\rm e}^{\frac{2\,a}{5}x}}.
\end{equation}

\begin{obs} Some comments:
\begin{enumerate}
\item These solutions are not new. We can find them in the article by Almendral and Sanjuán \cite{AlSa} along with a detailed discussion/description of the integrable (force free) cases, as well as an analysis of the physical behaviour of the Helmholtz oscillator with friction. In relation to the present article, the novelty consists in the method used to find the solutions. 
\item Another noteworthy element is the way in which the integrable cases are analysed. The entire analysis of the differential equations used in \cite{AlSa} to find the symmetries is here replaced by algebraic conditions that come naturally from the imposition of the existence of a Darboux integrating factor (which implies that $\sigma$ is a rational function).
\end{enumerate}
\end{obs}

\medskip

\subsection{Duffing-Van der Pol oscillator}
\label{dvpo}

A large amount of phenomena can be modelled by nonlinear oscillators and, among these, the Duffing-van der Pol oscillator occupies a prominent position, being able to describe from electrical circuits to low frequency oscillations of ion sound waves. The Duffing-van der Pol equation is a combination of the Duffing oscillator equation (an oscillator with a cubic non linearity) 
\begin{equation}
\ddot{x} +\alpha \dot{x} + {\omega_0}^2\,x + \beta\,x^3=0,
\label{2odeduff}
\end{equation}
that appears in numerous physics applications (see \cite{SaHeHe} and references therein) with the equation
\begin{equation}
\ddot{x} +\epsilon\,(1-x^2)\, \dot{x} + x =0,
\label{2odevdP}
\end{equation}
proposed by van der Pol to model oscillations in a vacuum tube triode circuit (see \cite{Pol,Pol2,Pol3}). When subjected to forced periodic excitation, these equations exhibit, in general, chaotic behavior (see \cite{HoRa}). So, the forced Duffing-Van der Pol equation is usually written as
\begin{equation}
\ddot{x} + \alpha \,(1-x^2)\,\dot{x}+\beta\,x+\gamma\,x^3=f\,\cos(\omega\, t).
\label{2odedvp}
\end{equation}
For applications in Plasma Physics \cite{KaNbOrTa,MiHiMoOr} the Duffing-van der Pol equation can be modified to
\begin{equation}
\ddot{x} + \epsilon \,(1-x^2)\,\dot{x}+x+\alpha\,x\,\dot{x}+\beta\,x^2+\gamma\,x^3=f\,\cos(\omega\, t),
\label{2odemdvp}
\end{equation}
where $\{\alpha,\,\beta,\,\gamma,\,\epsilon,\,f,\,\omega\}$ are parameters\footnote{The modified Duffing-van der Pol equation (\ref{2odemdvp}) is more general than (\ref{2odedvp}), so it can be used to find first integrals for the original Duffing-van der Pol equation equation as well.}. As the 2ODE (\ref{2odemdvp}) presents chaotic behaviour for arbitrary parameter values, it is of great importance for the study of the represented systems, the determination (if any) of sets of parameters values for which the 2ODE presents regular behaviour. In this way, the existence of a Lie symmetry, even of the type addressed in this work (i.e., non-local), guarantees the integrability of the system (for the particular set of values of the parameters that admits the symmetry). In order to check the existence (or not) of a Lie symmetry we can use the method developed in sections \ref{themet} and \ref{speci}. Since the procedure requires the 2ODE to be rational, we can make the transformation $\{{\rm e}^{i\,\omega \,t} \rightarrow x\}$, which leads to
\begin{equation}
 \left\{ \begin{array}{l}
\cos(\omega\, t) = \frac{1}{2}\,\left( x+ \frac{1}{x}\right),\\
x = y, \\
\dot{x} = i\,\omega \, y'\,x,  \\
\ddot{x} = -\omega^2x\,(y'+y''x),
\end{array} \right.
\end{equation}
and therefore, the 2ODE (\ref{2odemdvp}) is turned into        
\begin{equation}
y'' \!=\! {\frac {2i \epsilon \omega{x}^{2} y'(1\!-\!{y}^{2})\!+\!2ia\omega{x}^{2}yy' \!-\! 2\omega^2x^2z\!+\!2\gamma x{y}
^{3}\!+\!2\beta x\,{y}^{2}\!-\!2\omega {x}^{2}y'\!+\!2x\,y\!-\!f{x}^{2}\!-\!f}{2\omega^2{x}^{3}}}.
\label{2odemdvpr}
\end{equation}

Since the solutions of the set of undetermined coefficients  that leads to a symmetry (and, therefore, to an integrable case) are very numerous, in what follows we will present only two cases in which we can find Liouvillian first integrals for the forced regime ($f$ and $\omega$ are $\neq 0$) which, for arbitrary parameter values, is chaotic.

\medskip

\noindent
{\bf Case 1:} We find a solution for $\left\{ \alpha=0,\,\,\beta=0,\,\,\gamma=-\frac{1}{6}\,{e}^{2}+\frac{1}{6}\,ie\,k,\,\,\epsilon=\epsilon,\,\, f=f \right\}$, where $k \equiv \sqrt {4-{e}^{2}}$ (forced Duffing-van der Pol equation):
\begin{equation}
\label{sfdvdp}
\sigma= {\frac {i \left( 2\,e{y}^{2}+ik-e \right) }{2\omega x}}.
\end{equation}
The 2ODE becomes
\begin{equation}
\label{ode2mdvdpr}
y''=-{\frac { \left( 6\,i\epsilon\omega{x}^{2}{y}^{2}-6\,i\epsilon\omega{x}^{2}+6\,{\omega}^{2}{x}
^{2} \right) y'+i\epsilon kx{y}^{3}+{\epsilon}^{2}x{y}^{3}+3\,f{x}^{2}-6\,yx+3\,f}{6\,{x}
^{3}{\omega}^{2}}}.
\end{equation}
The symmetry is
\begin{equation}
{X_e}^{\!(1)}\! =\! {\rm e}^{\!\int_x\!\left[-{\frac {i \left( 2\,e{y}^{2}+ik-e \right) }{2\omega x}}\!\right]}\,\!\!\left(\!\partial_y\! -\! {\frac {i \left( 2\,e{y}^{2}+ik-e \right) }{2\omega x}}\,\partial_{y'}\!\right)\!.
\label{xcas2}
\end{equation}
Using the symmetry to determine the first integral and applying the inverse transformation
\begin{equation}
 \left\{ \begin{array}{l}
x = e^{i\omega\,t},\\
y = x, \\
y' = - \frac{i\,\dot{x}}{\omega \,e^{i\omega\,t}}, 
\end{array} \right.
\end{equation}
we obtain the first integral of the 2ODE:
\begin{eqnarray}
I&=& \left( \cos \left(kt/2 \right) +i\sin \left(kt/2 \right)  \right)  \left( i{\epsilon}^{2}k{x}^{3}-3\,ikx-3\,\epsilon{\omega}^{2}x+3\,i\cos \left( \omega t \right) fk+
 \right. \nonumber  \\ [3mm]
&& \left. 
3\,\cos \left( \omega t \right) \epsilon f+6\,\sin \left( \omega t \right) f\omega+2\,\epsilon{\omega}^{2}{x}^{3}-3\,\epsilon x-3(i\epsilon k+{\epsilon}^{2}+2{\omega}^{2}-2)\dot{x}+
\right. \nonumber  \\ [3mm]
&& \left. 
-2\,\epsilon{x}^{3}+{\epsilon}^{3}{x}^{3}+3\,ik{\omega}^{2}x \right) {{\rm e}^{\epsilon t/2}} / \left( ik\epsilon+{\epsilon}^{2}+2\,{\omega}^{2}-2 \right) 
\end{eqnarray}

\medskip

\noindent
{\bf Case 2:} We find a solution for $ \left\{ \alpha=\alpha,\beta=\beta,\gamma=-\frac {1}{3}-\frac {4}{3}\,{\frac {{\beta}^{2}}{{\alpha}^{2}}},\epsilon=\frac {1}{2}\,{\frac {\alpha}{\beta}}+2\,{\frac {\beta}{\alpha}},f=f \right\} $ (forced modified Duffing-van der Pol equation):
\begin{equation}
\label{sfdvdp2}
\sigma= {\frac {-i \left( 2\,{\alpha}^{2}\beta y-{\alpha}^{2}{y}^{2}-4\,{\beta}^{2}{y}^{2}+{\alpha}^{2} \right) }{2\alpha\beta\omega x}}.
\end{equation}
The 2ODE becomes
\begin{eqnarray}
y''&=& {\frac {iy' \left( 2\,{\alpha}^{2}\beta y-{\alpha}^{2}{y}^{2}-4\,{\beta}^{2}{y}^{2}+{\alpha}^{2}+4\,{\beta}^{2} \right) }{2wx\alpha\beta}}-\frac{y'}{x}+ \nonumber \\ [2mm]
&& {\frac {6\,\beta{y}^{2}x{\alpha}^{2}-2\,{\alpha}^{2}x{y}^{3}-8\,{\beta}^{2}x{y}^{3}-3\,f{x}^{2}{\alpha}^{2}+6\,yx{\alpha}^{2}-3\,f{\alpha}^{2}}{6{w}^{2}{\alpha}^{2}{x}^{3}}}.
\label{ode2mdvdpr2}
\end{eqnarray}
The symmetry is
\begin{equation}
{X_e}^{\!(1)}\! =\! {\rm e}^{\!\int_x\!\left[{\frac {i \left( 2\,{\alpha}^{2}\beta y-{\alpha}^{2}{y}^{2}-4\,{\beta}^{2}{y}^{2}+{\alpha}^{2} \right) }{2\alpha\beta\omega x}}\!\right]}\,\!\!\left(\!\partial_y\! +\!{\frac {i \left( 2\,{\alpha}^{2}\beta y-{\alpha}^{2}{y}^{2}-4\,{\beta}^{2}{y}^{2}+{\alpha}^{2} \right) }{2\alpha\beta\omega x}}\,\partial_{y'}\!\right)\!,
\label{xcas2}
\end{equation}
and, using it, we can determine the first integral that in the original coordinates can be written as
\begin{eqnarray}
I &=& \left( -\dot{x}+{\frac { \left({\alpha}^{2}+4\,{\beta}^{2}
 \right) {x}^{3}}{6\alpha \beta}}-{\frac {\alpha{x}^{2}}{2}}-{\frac {\alpha x}{2 \beta}} \right) {\rm e}^{\frac {2\beta}{\alpha}t}+  \nonumber  \\ [3mm]
&& \frac{f\alpha  \left( \sin \left( \omega t \right) \alpha \omega+2\,\cos \left( \omega t \right) \beta
 \right)}{{\alpha}^{2}{\omega}^{2}+4\,{\beta}^{2} } \, {\rm e}^{\frac {2\beta}{\alpha}t}.
\end{eqnarray}

\begin{obs} Some comments:
\begin{enumerate}
\item As mentioned, the forced (modified) Duffing-van der Pol oscillator presents chaotic behaviour for arbitrary values of the parameters and, thus, for the relations found between the parameters we have the regions where the behaviour of the system is regular.

\item As far as we know the first integrals found in cases 1 and 2 had not yet been discovered. This shows us that the method is very useful in determining first integrals.

\item Last but not least, the determination of regions of the parameters for which there is a symmetry (even if not local) implies that the system (in this region of parameters) is analytically integrable, that is, it presents an analytic first integral. This remains valid even if the first integral is not Liouvillian, i.e., even in cases where the 2ODEs have no Liouvillian first integrals.
\end{enumerate}
\end{obs}

\section{Conclusion}
\label{conclu}
\bigskip
\noindent

In this work we have introduced a new approach to search for symmetries (including the non-local ones) algorithmically for rational second order ordinary differential equations (rational 2ODEs). Our method is designed to deal with 2ODEs presenting at least one Liouvillian first differential invariant.

In very few words, our method is based in a few corner stones: We use a formal equivalence of the ``usual'' total derivative $d/dx$ and $D_x$ (the total derivative over the solutions).  

We have also restricted our method to the case where the integrating factor $\mu$ is given by: $$\mu = {\rm e}^{Z_0}\,\prod_i {p_i}^{n_i}$$ 
for the precise definition of the terms, please see section \ref{possalg}. 

This choice is not at all very restrictive, we have not actually found any case where it was not covered by this hypothesis. Furthermore, it is worth mentioning that the expression above is very general in what concerns elementary functions. 

Also, the above expression for $\mu$ has been proved to be the general expression for the integrating factor for the case of rational 1ODEs and, actually,  this very fact was what  motivated us to study it in the first place albeit, presently, for the case of rational 2ODEs. 

Another important characteristic of our approach is that it is designed to, basically, transform the core of the problem into algebraic calculations. This, for instance, apart turning all the calculations, as a general rule, much more efficient, makes the finding of regions of integrability, on the parameters present of the 2ODE, a simpler and fast endeavour (please, see section \ref{applic}).

It is worth to emphasise that, in principle, our approach is a semi-algorithm since a full algorithm must end in a finite number of steps. In our case, we can not determine the maximum degree of the polynomials $p$ and $q$ that form the $\sigma$ function. Nevertheless, up to the degree analysed, we can be sure that if we do not find the $\sigma$ function, then it does not exist within the range considered. 

{\bf Important remark:} Our (semi) algorithm becomes a (full) algorithm in the (very {\bf frequent}) cases pointed out (and demonstrated) in  section \ref{qeqn}, where $q | N$ or $q=u\,N$ where $u \in \{x,y,y'\}$.

Regarding the examples of applications hereby presented, we would like to point out the evidence for the efficiency of our approach embodied by the 2ODE on section \ref{dvpo}, where we believe we have found a new result.

We hope to have been able to demonstrate that our proposed algorithm is a valid contribution to the arsenal for dealing with rational 2ODEs. Apart from that, the mathematical results hereby presented can contribute to further the development of yet other lines of research, expanding the capabilities for integrating the class of 2ODEs under scrutiny here.


\end{document}